\documentclass[12pt,a4paper]{article}
\usepackage{graphicx} 
\usepackage{epsfig}
\usepackage{amssymb}
\usepackage{amsmath}
\usepackage{amsfonts}
\usepackage{mathrsfs}
\usepackage[english]{babel}   
\usepackage{fancyhdr}
\usepackage{lastpage}
\usepackage[normal]{subfigure}
\usepackage{color}
\usepackage[dvipsnames]{xcolor}
\usepackage{inputenc}
\usepackage{dashrule}
\usepackage{float}
\usepackage{soul}

\newtheorem{comentario}{{\bf Note}}[section]

\newtheorem{remark}{Remark}[section]

\def\QH{  { {\mathcal U}_{\mathcal Z}^H}}

\def\bA{{\bf A}}
\def\bH{{\bf H}}
\def\bG{{\bf G}}
\def\bF{{\bf F}}
\def\bU{{\bf U}}

\def\vu{\boldsymbol{u}}

\def\dxi{\Delta x_{i}}
\def\dx{\Delta x}
\def\dt{\Delta t}
\def\dxip{\Delta x_{i+\frac12}}
\def\dxim{\Delta x_{i-\frac12}}

\def\a{\alpha}
\def\b{\beta}
\def\r{\rho}
\def\O{\Omega}
\def\D{\Delta}
\def\g{\nabla}
\def\p{\partial}
\def\wtd{\widetilde}
\def\dsum{\displaystyle\sum}

\newcommand{\abs}[1]{\left\vert #1 \right\vert}

\definecolor{deepskyblue4}{rgb}{0,0.41,0.55}
\definecolor{dodgerblue4}{rgb}{0.06,0.31,0.55}
\definecolor{cadetblue2}{rgb}{0.56,0.9,0.93}
\definecolor{paleturquoise2}{rgb}{0.68,0.93,0.93}
\definecolor{deepskyblue3}{rgb}{0,0.6,0.8}
\definecolor{deepskyblue2}{rgb}{0,0.7,0.93}
\definecolor{royalblue3}{rgb}{0.23,0.37,0.8}
\definecolor{chartreuse3}{rgb}{0.4,0.8,0}
\definecolor{green3}{rgb}{0,0.8,0}
\definecolor{cyan3}{rgb}{0,0.8,0.8}
\definecolor{orchid}{rgb}{0.85,0.44,0.84}
\definecolor{orange}{rgb}{1,0.65,0}
\definecolor{darkorange}{rgb}{1,0.55,0}
\definecolor{sienna1}{rgb}{1,0.51,0.28}
\definecolor{red2}{rgb}{0.93,0,0}
\definecolor{red1}{rgb}{1,0,0}
\definecolor{firebrick2}{rgb}{0.93,0.17,0.17}
\definecolor{yellow2}{rgb}{0.93,0.93,0}
\definecolor{gold1}{rgb}{1,0.84,0}
\definecolor{orchid}{rgb}{0.85,0.44,0.84}

\title{Multilayer shallow water models with \\ locally variable number of layers\\ 
and semi-implicit  time discretization}

\author{Luca Bonaventura $^{(1)}$,  Enrique D. Fern\'andez-Nieto $^{(2)}$,\\ Jos\'e Garres-D\'iaz $^{(2)}$ and Gladys Narbona-Reina $^{(2)}$}

\begin{document}
\maketitle

\begin{center} 
{\small
$^{(1)}$ MOX -- Modelling and Scientific Computing, \\
Dipartimento di Matematica, Politecnico di Milano \\
Via Bonardi 9, 20133 Milano, Italy\\
{\tt luca.bonaventura@polimi.it}
}
\end{center}

\begin{center}
{\small $^{(2)}$
 IMUS \& Departamento de Matem\'atica Aplicada I. \\
 ETS Arquitectura, Universidad de Sevilla\\
Avda. Reina Mercedes 2, 41012, Sevilla, Spain\\
{\tt edofer@us.es, jgarres@us.es, gnarbona@us.es}
}
\end{center}
\date{}


\vspace*{0.5cm}


\vspace*{0.5cm}

\pagebreak

\begin{abstract}
We propose an extension of the discretization approaches for multilayer shallow water models, aimed at making them more flexible and efficient for realistic applications to coastal flows. A novel discretization approach is proposed, in which the number of vertical  layers and their distribution
are allowed to change in different regions of the computational domain. Furthermore, semi-implicit  schemes are employed for  the time discretization, leading to a significant efficiency improvement for subcritical regimes.  We show that, in the typical regimes in which the application of multilayer shallow water models is justified, the resulting discretization does not introduce any major spurious feature and allows again to reduce substantially the computational cost in areas with complex bathymetry. As an example of the potential of the proposed technique,
an application to a sediment transport problem is presented, showing a remarkable improvement with respect to standard discretization approaches.  
\end{abstract}

\section{Introduction}
\label{intro}

\indent

Multilayer shallow water models have been first proposed in \cite{audusse:2005} to account 
for the vertical structure in  the simulation of large scale geophysical flows. They have been later
extended and applied in \cite{audusse:2008}, \cite{audusse:2011}, \cite{audusse:2011b}.  This multilayer model was applied in \cite{audusse:2011c} to study movable beds by adding an Exner equation. A different formulation, to which we will refer in this paper, was proposed in \cite{fernandez:2014}, which has several peculiarities with respect to   previous multilayer models.
The model proposed in \cite{fernandez:2014}  is derived from the weak form of the full Navier-Stokes system, by assuming a discontinuous profile of velocity, and the solution is obtained as a particular weak solution of the full Navier-Stokes system. The vertical velocity is computed in a postprocessing step based on the incompressibility condition, but accounting also  for the mass transfer terms between the internal layers. In \cite{fernandez:2016}, this multilayer approach is applied to dry granular flows, for which an accurate approximation of the  vertical flow structure is essential  to approximate the velocity-pressure dependent viscosity.
  
Multilayer shallow water models can be seen as an alternative
to more standard approaches for vertical discretizations, such as natural height coordinates,
  (also known as $z-$coordinates in the literature on numerical modelling of atmospheric and oceanic flows), employed e.g. in  \cite{bonaventura:2000}, \cite{bryan:1969},    \cite{casulli:1994}, terrain following coordinates (also known as $\sigma-$coordinates in the literature), see e.g. \cite{haidvogel:1991},
and isopycnal coordinates, see e.g. \cite{bleck:1986}, \cite{casulli:1997}. Each technique has its own advantages and shortcomings,
as highlighted in  the discussions and reviews in \cite{adcroft:1997}, \cite{bonaventura:2000}, \cite{bonaventura:2012}, \cite{haney:1991}.
Multilayer approaches are appealing, because they  share some of the advantages of $z-$coordinates, such as the absence of metric terms in the model equations, while not requiring special treatment of the lower boundary. On the other hand, multilayer approaches share one of the main disadvantages of $\sigma-$coordinates, since they require, at least in the formulations employed so far, to use the same number of layers independently   of the fluid  depth. 
Furthermore, an implicit regularity assumption on the lower boundary is required, in order to avoid that
too steeply inclined layers arise, which would contradict the fundamental hydrostatic assumption underlying the model.

In this work, we propose two concurrent strategies to make multilayer models more efficient and fully competitive with their $z-$ and $\sigma-$coordinates counterparts. On one hand, we propose a novel discretization approach,
 in which the  number of vertical layers can vary over the computational domain. We show that, in the typical regimes in which the application of multilayer shallow water models is justified, the resulting discretization does not introduce significant errors and allows to reduce substantially the computational cost in areas with complex bathymetry. Thus making multilayer approach fully competitive with $z-$coordinate discretizations for large scale, hydrostatic flows. Furthermore, efficient semi-implicit discretizations are applied for the first time to this kind of models, allowing to achieve the same kind of computational gains that have been obtained for other vertical discretization approaches. In this paper, for simplicity, we have restricted our attention to constant density flows. An extension to variable density problems in the Boussinesq regime  will be presented in a forthcoming paper.
 
 In section \ref{se:model}, the equations defining the multilayer shallow water models of interest will be reviewed. In section \ref{spatial}, the  spatial discretization is introduced in a simplified framework, showing how the number of layers can be allowed to vary over the computational domain. In section \ref{si_td}, some semi-implicit time discretizations are introduced for the model with a variable number of layers. Results of a number of numerical experiments are reported in section \ref{tests}, showing the significant efficiency gains that can be achieved by combination of these two techniques. Some conclusions and perspectives for future work are presented in section \ref{conclu}. 

\bigskip
\section{Multilayer shallow water models}
\label{se:model}
\begin{figure}[!h]
\begin{center}
\includegraphics[width=1\textwidth]{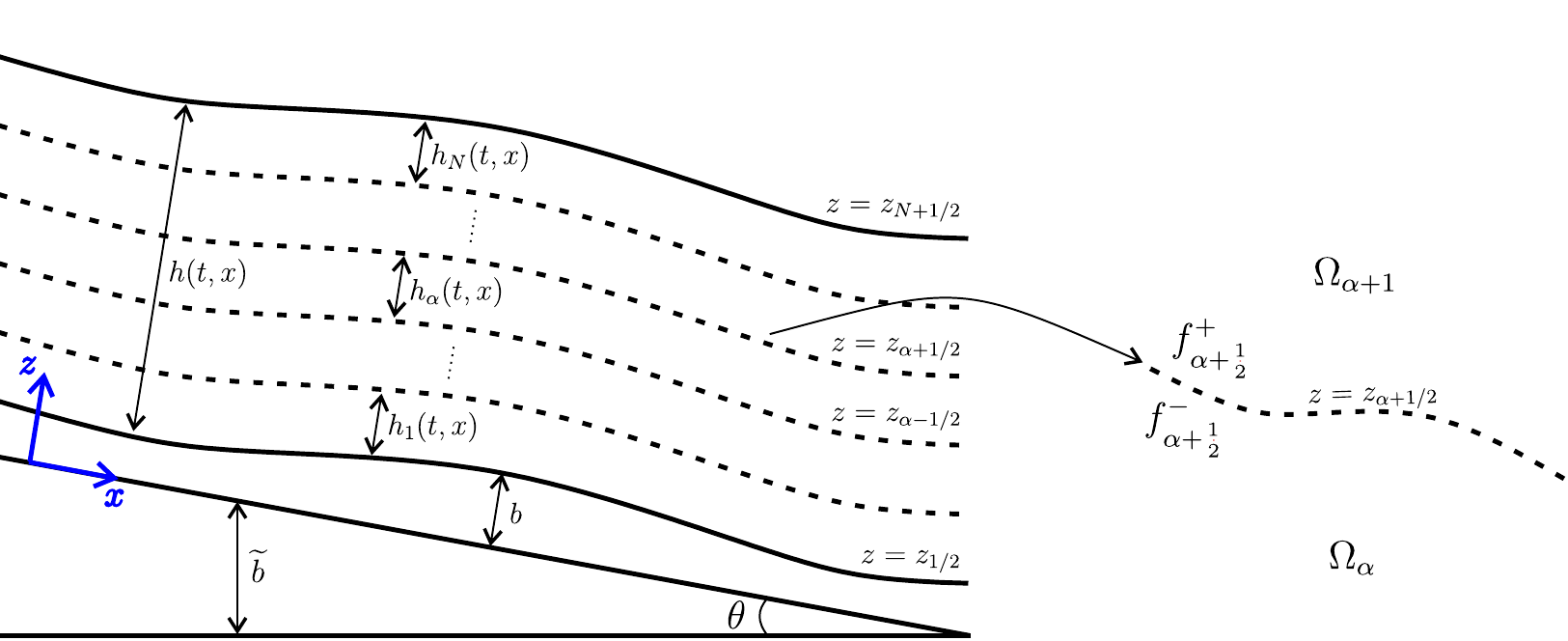}
\end{center}
 \caption{\label{fig:Multilayers} \it{Sketch of the domain and of its subdivision in a constant number of   layers.}}
 \end{figure}

We consider the multilayer shallow water model described pictorially in Figure \ref{fig:Multilayers}. In this approach,
  $N$ subdivisions $\O_\a,\ \a=1,\dots,N$ of the domain $\O$ are introduced 
  in the vertical direction. We denote by $h_\a$ the height of the layer $\a$ and by $h = \sum_{\a=1}^{N}h_\a$ the total height. Note that $\O = \bigcup^N_{\a=1} \O_\a$ and that each subdomain $\O_\a$ is delimited by 
  time dependent interfaces $\Gamma_{\a\pm\frac12}(t),$ that are assumed to be represented by the one valued functions $z = z_{\a\pm\frac12}(x,y,t)$.  For a function $\smash f$ and for $\smash\a=0,1,...,N$, we also define, as in
 \cite{fernandez:2014}, 
  $$ f_{\a+\frac{1}{2}}^{-} := ( f_{|_{\O_{\a}(t)}})_{|_{\Gamma_{\a+\frac{1}{2}}(t)}} \text{\;\;and\;\;}  f_{\a+\frac{1}{2}}^{+} := (f_{|_{\O_{\a+1}(t)}})_{|_{\Gamma_{\a+\frac{1}{2}}(t)}}.$$
Obviously, if the function $\smash f$ is continuous,
  $$ f_{\a+\frac{1}{2}} := f_{|_{\Gamma_{\a+\frac{1}{2}}(t)}} = f_{\a+\frac{1}{2}}^+ = f_{\a+\frac{1}{2}}^-.$$
  
\noindent Following \cite{fernandez:2016}, \cite{fernandez:2014},  the equations describing this multilayer approach can be  written 
 for $\smash\a = 1,\dots,N$ as
\begin{align}
\label{eq:FinalModel}
\begin{array}{l}
\p_{t}h_{\a} + \g\!_{x}\cdot(h_{\a}\vu_{\a}) = G_{\a+\frac{1}{2}} - G_{\a-\frac{1}{2}}, \\
\\
\p_{t}\left(h_{\a}\vu_{\a}\right) \;+\; \g\!_{x}\cdot\left( h_{\a}\vu_{\a}\otimes\vu_{\a}\right) \, + \\[4mm]
\quad+ gh_{\a}\g_{x} \left(b+h\right) \, =\, \dfrac{1}{\r_0}\left(\boldsymbol{K}_{\a-\frac{1}{2}} - \boldsymbol{K}_{\a+\frac{1}{2}}\right)\; +\\[4mm]
\quad+\;\dfrac{1}{2} G_{\a+\frac{1}{2}}\left(\vu_{\a+1} + \vu_{\a}\right) \;-\; \dfrac{1}{2} G_{\a-\frac{1}{2}}\left(\vu_{\a} + \vu_{\a-1}\right).
\\
\end{array}
\end{align}
\\
Here, we consider a fluid with constant density $\rho_0,$
  $(\vu_\a, w_\a) \in \mathbb{R}^3$ is the velocity in the layer $\a,$  $b=b(x,y)$ is a function describing the bathymetry (which is assumed to be constant in time, so that  $\p_t b = 0 $) 
and the terms $ \boldsymbol{K}_{\a+\frac{1}{2}} $ model the shear stresses between the layers.  Notice that the atmospheric pressure has assumed to be zero.
The vertical velocity profile  is recovered from the integrated incompressibility condition, obtaining for $\smash\a=1,...,N$ and $z \in (z_{\a-\frac{1}{2}}, z_{\a+\frac{1}{2}})$,
\begin{equation}
\label{eq:vert_vel}
 w_{\a}(t,x,z) = w_{\a-\frac{1}{2}}^+(t,x) \;-\; (z-z_{\a-\frac{1}{2}})\g\!_{x}\cdot \vu_{\a}(t,x),
  \end{equation}
where
$$
  w_{\a+\frac{1}{2}}^{+} = (\vu_{\a+1} - \vu_{\a})\;\g\!_{x}z_{\a+\frac{1}{2}} \;+\; w_{\a+\frac{1}{2}}^{-},
$$
and
$$  
w_{\a+\frac{1}{2}}^{-} =w_{\a-\frac{1}{2}}^{+} - h_{\a}\g\!_{x}\cdot\vu_{\a} 
\quad \mbox{with} \quad 
w_{\frac12}^+=u_1 \g_x b - G_{\frac12}.
$$
Since we are focusing in this work mostly on subcritical flows, there is no 
special reason to choose discharge rather than velocity as a model variable. Therefore, we rewrite the previous 
system as 
\begin{align}
\label{eq:FinalModel_u}
\begin{array}{l}
\p_{t}h_{\a} + \g\!_{x}\cdot(h_{\a}\vu_{\a}) = G_{\a+\frac{1}{2}} - G_{\a-\frac{1}{2}}, \\[4mm]
h_{\a}\p_{t} \vu_{\a} \;+\; h_{\a}\vu_{\a}\cdot \nabla \vu_{\a}  \,  + gh_{\a}\g_{x} \left(b+h\right) = \\[4mm]
\quad  \quad \quad \quad =\dfrac{1}{\r_0}\left(\boldsymbol{K}_{\a-\frac{1}{2}} - \boldsymbol{K}_{\a+\frac{1}{2}}\right)\; 
+\left ( G_{\a+\frac{1}{2}}\Delta \tilde{\vu}_{\a+\frac 12} +G_{\a-\frac{1}{2}}\Delta \tilde{\vu}_{\a-\frac 12}\right),\\
\end{array}
\end{align}
where  $\Delta \tilde{\vu}_{\a+\frac 12}=( {\vu}_{\a+1}-{\vu}_{\a})/2.$
From the derivation in \cite{fernandez:2014}, it follows that\\
\begin{eqnarray}
\label{eq:FinalModel2}
G_{\a+\frac{1}{2}}& = &\p_{t}z_{\a+\frac{1}{2}} + \vu_{\a+1}\cdot\g\!_{x}z_{\a+\frac{1}{2}} - w_{\a+\frac{1}{2}}^{+} = \p_{t}z_{\a+\frac{1}{2}} + \vu_{\a}\cdot\g\!_{x}z_{\a+\frac{1}{2}} - w_{\a+\frac{1}{2}}^{-},\nonumber 
\\
\boldsymbol{K}_{\a+\frac{1}{2}} &= &- \mu_{\a+\frac{1}{2}} \QH_{\a+\frac12} ,
\end{eqnarray}
where $\mu$ denotes the dynamic viscosity   and $\QH_{\a+\frac12}$ is an approximation of $\p_z u_{\a}$ at $\Gamma_{\a+\frac12}$.
We then define the vertical partition of the domain, setting $h_{\a} = l_{\a}\,h \,$ where, for $ \a = 1,\cdots,N, $ $l_{\a}\ $ are positive constants
such that
\begin{equation*}
 \label{eq:layersum}
\dsum_{\a=1}^N l_{\a} = 1.
\end{equation*}

\noindent Note that model (\ref{eq:FinalModel})  consists of  $2N$ equations in the unknowns 
$$ h, \{ \vu_{\a}\}_{\a=1,\dots,N},  \ \ \ \ \{G_{\a+\frac{1}{2}} \}_{\a=1, \dots,N-1}. $$
However, the mass transfer terms can be rewritten as  
$$
G_{\a+\frac{1}{2}} = \p_{t}z_{\a+\frac{1}{2}} + \frac{\vu_{\a}+ \vu_{\a+1}}{2}\p\!_{x}z_{\a+\frac{1}{2}} - w_{\a+\frac{1}{2}},
\quad \mbox{where} \quad w_{\a+\frac{1}{2}}= \frac{w_{\a+\frac{1}{2}}^+ + w_{\a+\frac{1}{2}}^-}{2}.
$$
As a consequence, the system has $2N$ unknowns, now corresponding to the total height $h$, the velocity  $ \{\vu_{\a} \}_{\a=1,\dots,N}$
in each layer and the averaged vertical velocity at each internal interface $\{ w_{\a+\frac{1}{2}} \}_{\a=1,\dots, N-1}$.
By combining the continuity equations, the system can be rewritten with $N+1$ equations and unknowns. The unknowns of the reduced system are the total height $h $ and the velocity in each layer, $\vu_{\a},$ for 
$\a=1,\dots,N.$
\noindent Note that $G_{\a+\frac{1}{2}}$ can be written, by summing the mass equations from 1 to $\a$, as
\begin{align}
\label{eq:cont_mass_al}
G_{\a+\frac{1}{2}} = G_{\frac{1}{2}} + \dsum_{\b=1}^\a \left(\p_{t}h_{\b} + \g_{x}\cdot( h_{\b}\vu_{\b})\right).
\end{align}

\noindent Moreover, for the special case $\a=N$ and   $G_{1/2} = G_{N+\frac{1}{2}}=0$, the above equation leads to
\begin{align*}
\label{eq:cont_mass}
\p_{t}h  + \g_{x}\cdot\Biggl(  h\dsum_{\b=1}^N l_{\b}\vu_{\b}\Biggr) = 0.
\end{align*}

\noindent By introducing this in the mass equation we obtain
\begin{equation}\label{eq:cont_mass_G}
 G_{\a+\frac{1}{2}} = \  \dsum_{\b=1}^\a \Biggl( \g_{x}\cdot\left(hl_{\b}\vu_{\b}\right) -  l_{\b}\dsum_{\gamma=1}^N\g_{x}\cdot\left(l_{\gamma}h \vu_{\gamma}\right)\Biggr).
\end{equation}

\noindent
  Assuming  also   $\p_t b = 0, $  
 system \eqref{eq:FinalModel_u}-\eqref{eq:FinalModel2} is finally re-written as
\begin{eqnarray}
\label{eq:sistema3D_u}
 \p_{t}\eta  &+& \g_{x}\cdot\Biggl(  h\dsum_{\b=1}^N l_{\b}\vu_{\b}\Biggr) = 0, \nonumber  \\[2mm]
 \p_{t} \vu_{\a} &+&  \vu_{\a}\cdot \nabla \vu_{\a}  \,  + g \g_{x} \eta  =  \\[4mm]
 && \qquad =\dfrac{\boldsymbol{K}_{\a-\frac{1}{2}} - \boldsymbol{K}_{\a+\frac{1}{2}}}{\r_0 h_{\a}}\; 
+ \dfrac{ G_{\a+\frac{1}{2}}\Delta \tilde{\vu}_{\a+\frac 12} +G_{\a-\frac{1}{2}}\Delta \tilde{\vu}_{\a-\frac 12}}{h_{\a}}, \nonumber 
\end{eqnarray}
for $\a=1,\cdots,N.$ Here, we have set as customary in the literature $\eta = b + h $. The transport equation for a passive scalar can be coupled to the previous continuity and momentum equation, in such a way as to guarantee compatibility with the  continuity equation in the sense of
\cite{gross:2002}. 
  If $\rho_{\a}$ denotes the average density of the passive scalar in $\Omega_{\a}$, it verifies the following tracer equation:
$$
\p_{t}\left(\r_\a h_{\a}\right) + \g\!_{x}\cdot(\r_\a h_{\a}\vu_{\a}) = \r_{\a+1/2}G_{\a+\frac{1}{2}} - \r_{\a-1/2}G_{\a-\frac{1}{2}},  
$$
where
$$
 \r_{\a+1/2}=\frac{\r_{\a}+\r_{\a+1}}{2} + \frac{1}{2} \mbox{sgn}(G_{\a+\frac{1}{2}})(\r_{\a+1}-\r_{\a}).
$$

In principle, 
any appropriate turbulence and friction model can   be considered to define the turbulent fluxes $\boldsymbol{K}_{\a+\frac{1}{2}}$, $\a=0,\dots,N$.  Here, we have employed a
parabolic turbulent viscosity profile and friction coefficients derived from a logarithmic wall law:
$$
\nu = \dfrac{\mu}{\r_0} = \kappa\,u^*\,(z-b)\left(1-\dfrac{z-b}{h}\right),
$$
where $\kappa=0.41$ is the von Karman constant, $u^*=\sqrt{\tau_b/\rho}$ is the friction velocity and  $\tau_b $ denotes the shear stress.    In order to approximate this turbulence model we set for $\a=1, \dots, N-1$:
$$
\boldsymbol{K}_{\a+\frac{1}{2}}= \mu_{\a+\frac12} \frac{u_{\a+1}-u_{\a}}{(h_{\a}+h_{\a+1})/2}, \quad\text{with}\quad \mu_{\a+\frac12} = \r_0 \kappa \, u^*_{\a+\frac{1}{2}} \, \left( \sum_{\beta=1}^{\a} l_{\beta}h \right) \left( \sum_{\gamma=\a+1}^{N} l_{\gamma} \right).
$$
Trivially, $\nu_{\a+\frac12} = \mu_{\a+\frac12}/\r_0$. For $\a=0$ and $\a=N$, standard quadratic models for bottom  and wind stress are considered. We then set 
$$
\boldsymbol{K}_{1/2} =  -C_{f}  \|\vu_1\|\vu_1, \quad \boldsymbol{K}_{N+1/2}  = -C_{w} \|\vu_{w}-\vu_N\| (\vu_{w}-\vu_N),
$$
 where $\vu_w$ denotes the wind velocity and $C_w$ the friction coefficient between at the free surface.  The friction coefficient $C_f $ is defined, according to the derivation in \cite{decoene:2009}, as:
 \begin{equation} \label{eq_def_cf}
C_f = \  \kappa^2\dfrac{\left(1-\dfrac{\D z_r}{h}\right)}{\left( \ln{\left(\dfrac{\D z_r}{\D z_0} \right)} \right)^2},
\end{equation}
where $\D z_0$ is the roughness length and $\Delta z_r$ is the length scale for the bottom layer.   Under the assumption that $\Delta z_0<< \Delta z_r $ it can be seen  that
$$
\frac{ u_t}{u^*}\approx  \frac{1}{\kappa} \ln {\big (\frac{z-b}{\Delta z_0} \big )}
 $$
where $u_t$ is the  tangential velocity.  In practice, we identify $ \displaystyle u_t$ with $u_1$, the horizontal velocity of the  layer closest to the bottom, in the multilayer model. The definition of $C_f$ given by equation (\ref{eq_def_cf}) is deduced by using previous  relation of the ratio between $u_1$ and $u^*$ (see \cite{decoene:2009}). Then, we set
$$
u^*_{\a+\frac{1}{2}}= \frac{u_1 \kappa}{ \displaystyle \ln {\big ( \sum_{\beta=1}^{\a} l_{\beta} h /  \Delta z_0 \big )}},
$$
in the definition of $\boldsymbol{K}_{\a+\frac{1}{2}}$.

\section{Spatial discretization with variable number of layers}
\label{spatial}

The multilayer shallow water model \eqref{eq:sistema3D_u} can be discretized in principle with any spatial discretization approach. 
For simplicity,   we present the proposed discretization approach in the framework of
 simple finite volume/finite difference discretization on a 
staggered Cartesian mesh with C-grid staggering. A discussion of the advantages
of this approach for large scale geophysical models can be found in
 \cite{durran:2013}. The C-grid staggering also has the side benefit of providing a more compact
structure for the system of equations that is obtained when a semi-implicit method is applied for time discretization. In order to further simplify the presentation, we only introduce the discretization for 
an $x-z$ vertical slice, even though any of the
methods presented in the following can be easily generalized to the full three dimensional case.
Generalization to structured and unstructured meshes can be obtained e.g. by the approaches proposed in \cite{casulli:1994} and \cite{bonaventura:2005}, \cite{casulli:2000},  \cite{casulli:2002}, respectively, but higher order  methods  such as those of \cite{tumolo:2015}, \cite{tumolo:2013} could also be applied.
It is to be remarked that the choice of a staggered mesh is by no means necessary and that the approach
proposed below to handle a variable number of layers can be easily extended to colocated meshes as well.
 
The   solution domain will then coincide with an interval $[0,L], $
that is assumed to be subdivided into control volumes $V_i, i=1,\dots,M.$ 
The step in the mesh is defined by $\dxi = x_{i+\frac12} - x_{i-\frac12}$ and $\dxip = x_{i+1}- x_{i}$, where $x_{i+\frac12} = (x_i + x_{i+1})/2$ as usual.  The discrete free surface    variables $\eta_i$ are defined  at integer locations corresponding to  the centers of the control volumes, 
while the discrete velocities $u_{\a,i+\frac12}$  are defined at the half-integer locations $i+1/2.$
As suggested in \cite{gross:2002}, the value of $h_{i+\frac12}$ is  taken to be that of the control volume located upwind of the volume edge. 
 
 The vertical number of layers employed, in the approach proposed in  \cite{fernandez:2014}, is a discretization parameter whose choice depends on the desired accuracy in the approximation of
the vertical structure of the flow. In order to make this type of model more flexible and more efficient,
we propose to allow for a number of vertical layers that is not constant throughout the domain. The transition between  regions with different numbers of layers is assumed to take place at the center of a control volume $V_i,$  so that one may have different $N_{i+\frac12}$ for $i=0,\dots,M $ and as a consequence, the discrete layer thickness coefficients $l_{\a,i+\frac12}$ are also defined
at  the half-integer locations $i+1/2.$
The number of layers considered at the cell center for the purpose of the discretization of the
tracer equation are defined as $N_i=\max{\{N_{i-\frac12},N_{i+\frac12}\}} $ and the discrete layer thickness coefficients at integer locations $l_{\a,i}  $ are taken to be equal to those at the neighbouring half-integer
location with larger number of layers.
 We will also assume that, 
whenever for some $i+\frac12$ one has, without loss of generality,
$N_{i-\frac12} > N_{i+\frac12}, $ then for any $\b = 1,\dots, N_{i+\frac12} $ 
there exist 
\begin{equation}
\label{interf}
  1\leq \a^-_{i-\frac12}(\b) \leq \a^+_{i-\frac12}(\b) \leq N_{i-\frac12} \ \ \ {\rm such \  that}  \ \ \ l_{\b,i+\frac12}=\sum_{\a= \a^-_{i-\frac12}(\b)  }^{\a^+_{i-\frac12}(\b)}  l_{\a,i-\frac12}.
  \end{equation}
This allows a more straightforward implementation of the numerical approximation of   horizontal advection in the velocity and in the tracer equation, which are the only ones involving a horizontal stencil.
Finally, again for simplicity of the implementation and without great loss of generality,
it is assumed that if $N_{i-\frac12} \neq N_{i+\frac12} $ one has 
 $N_{i-\frac32}  = N_{i-\frac12} $ as well as  $N_{i+\frac32}  = N_{i+\frac12}.$
 
 \begin{figure}[!h]
\begin{center}
\includegraphics[width=0.5\textwidth]{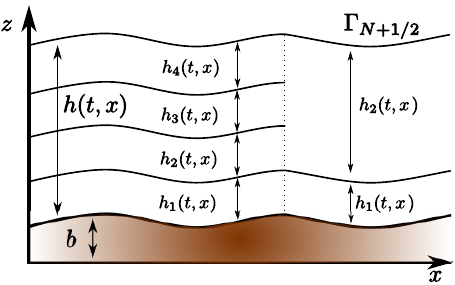}
\end{center}
 \caption{\label{fig:Multilayers_part} \it{Sketch of the domain and of its subdivision in a variable number of   layers.}}
 \end{figure}
 A sample configuration of this kind is depicted in figure \ref{fig:Multilayers_part}. 
  Notice that a dependence of the number of layers on time
    could also be introduced, in order to adapt the global maximum number of layers to
 the flow conditions, but this has not been done in the present implementation.

  \section{Semi-implicit time discretizations}
 \label{si_td} \indent

The previous definitions yield a space discretization that can be easily coupled to any
time discretization that yields a stable fully discrete space time scheme. For example,
a time discretization by a third order Runge Kutta scheme has been employed as a reference in
the numerical tests presented in section \ref{tests}. However, we will focus here on
semi-implicit time discretization approaches aimed at reducing the computational
cost in subcritical regime simulations.

With this goal, it is immediate to notice that the formal structure of system \eqref{eq:sistema3D_u}
 is entirely analogous to that of the three dimensional hydrostatic system
considered in  \cite{casulli:1994}, \cite{casulli:1992},  so that we can build  semi-implicit time discretizations along the same lines, i.e. by treating implicitly the velocity in the continuity equation and the free surface gradient in the momentum equation. In the following, we present first a more conventional   time discretization based
on the off-centered trapezoidal rule (or $\theta$-method, see e.g. \cite{lambert:1991}) and then a more
advanced Implicit-Explicit Additive Runge Kutta second order method (IMEX-ARK2). 

\subsection{A $\theta$-method time discretization}
Following   \cite{casulli:1994}, we first consider a semi-implicit discretization based on  the $\theta$-method, 
which can
be defined for  a generic ODE system $\mathbf{y}' = \mathbf{f}(\mathbf{y},t) $ as 
$$
\mathbf{y}_{n+1} = \mathbf{y}_n + \Delta t
\,\left[\theta\mathbf{f}(\mathbf{y}_{n+1},t_{n+1}) + (1-\theta)\mathbf{f}(\mathbf{y}_n,t_n)\right],
$$
where  $ \Delta t $ denotes the time step and $\theta \in [0,1]$ is a implicitness parameter. 
 If $\theta \geq 1/2$ the method is unconditionally stable and the numerical diffusion introduced by the method increases when increasing $\theta$. For $\theta =1/2$  the second order Crank-Nicolson method 
 is obtained.  In practical applications, $\theta$ is usually chosen just slightly larger than $1/2, $
 in order to allow for some damping of the fastest linear modes and nonlinear effects. We then proceed to describe the time discretization of system \eqref{eq:sistema3D_u} based on the $\theta$-method. 
  
  For control volume $i,$ the continuity equation  in \eqref{eq:sistema3D_u} is then discretized
 as 
 
 \begin{eqnarray}\label{eq:mass_discrete_theta}
&&\dxi \eta_{i}^{n+1}+ \,\theta\dt \left( \dsum_{\b=1}^{N_{i+\frac12}}
 l_{\b,i+\frac12} h^n_{i+\frac12} u_{\b,i+\frac12}^{n+1} \, - \, 
  \dsum_{\b=1}^{N_{i-\frac12}} l_{\b,i-\frac12} h^n_{i-\frac12}u_{\b,i-\frac12}^{n+1}\right)  \\
&&=   \dxi\eta_{i}^{n} -\,\left(1-\theta\right)\dt \left( \dsum_{\b=1}^{N_{i+\frac12}} l_{\b,i+\frac12} h^n_{i+\frac12} u_{\b,i+\frac12}^{n} \, - \,  \dsum_{\b=1}^{N_{i-\frac12}} l_{\b,i-\frac12} h^n_{i-\frac12}u_{\b,i-\frac12}^{n}\right).
\nonumber
\end{eqnarray}
 It can be noticed that the dependency on $h$ has been frozen at time level $n$ in order to
 avoid solving a nonlinear system at each timestep. As shown in \cite{casulli:1994}, \cite{tumolo:2015},
 this does not degrade the accuracy of the method, even in the case of a full second order discretization
 is employed.
 For nodes $i+\frac 12,$ the momentum equations for $\a=2,...,N_{i+\frac12}-1$ in 
 \eqref{eq:sistema3D_u} are then discretized as
 
 \begin{eqnarray}
\label{eq:vel_discrete_theta}
&&u_{\a,i+\frac12}^{n+1} + g\theta\frac{\dt} {\dxip}(\eta_{i+1}^{n+1} - \eta_{i}^{n+1}) \\
&&-\frac{\dt \theta}
{l_{\a,i+\frac12}\,h_{i+\frac12}^{n}}
\left(\nu_{\a+\frac12,i+\frac12}^{n}\,\dfrac{u_{\a+1,i+\frac12}^{n+1}-u_{\a,i+\frac12}^{n+1}}{l_{\a+\frac12,i+\frac12}h_{i+\frac12}^{n}} - \nu_{\a-\frac12,i+\frac12}^{n}\,\dfrac{u_{\a,i+\frac12}^{n+1}-u_{\a-1,i+\frac12}^{n+1}}{l_{\a-\frac12,i+\frac12}h_{i+\frac12}^{n}} \right) \nonumber\\
&&= u_{\a,i+\frac12}^{n}
 + \dt  {\cal A}_{\a,i+\frac 12}^{u,n}   -   g(1-\theta)\frac{\dt}{\dxip}\left(\eta_{i+1}^{n} - \eta_{i}^{n}\right) \nonumber
 \\ 
&&+\frac{\dt (1-\theta)}{l_{\a,i+\frac12}\,h_{i+\frac12}^{n}}\left(\nu_{\a+\frac12,i+\frac12}^{n}\,\dfrac{u_{\a+1,i+\frac12}^{n}-u_{\a,i+\frac12}^{n}}{l_{\a+\frac12,i+\frac12}h_{i+\frac12}^{n}} 
- \nu_{\a-\frac12,i+\frac12}^{n}\,\dfrac{u_{\a,i+\frac12}^{n}-u_{\a-1,i+\frac12}^{n}}{l_{\a-\frac12,i+\frac12}h_{i+\frac12}^{n}} \right) \nonumber\\ 
 &&+\frac{\dt} {\dxip l_{\a,i+\frac12}\,h_{i+\frac12}^{n}}\left(\Delta \wtd{u}_{\a+\frac12,i+\frac12}^{\,n}
  \mathcal{G}_{\a+\frac12,i+\frac12}^{n} 
+\Delta \wtd{u}_{\a-\frac12,i+\frac12}^{\,n} \mathcal{G}_{\a-\frac12,i+\frac12}^{n}\right), \nonumber
\end{eqnarray}
where  $\Delta \wtd{u}_{\a+\frac12,i+\frac12}^{\,n} =(  u_{\a+1,i+\frac12}^{\,n} -u_{\a,i+\frac12}^{\,n})/2, $
$\mathcal{G}_{\a+\frac12,i+\frac12}^{n} $ denotes a discretization of the mass transfer term and
${\cal A}_{\a,i+\frac 12}^{u,n} $ denotes some spatial discretization of the velocity advection term.
In the present implementation, an upstream based second order scheme is employed for this term.
 Notice that, to define this advection term, velocity values from different layers may have to be employed,
 if some of the neighbouring volumes has a number of layers different from that at $i+\frac 12.$ For
 example, assuming again without loss of generality $N_{i-\frac12} > N_{i+\frac12} $  and $u_{\b,i+\frac12}^{n} > 0 $ and using the notation in \eqref{interf},
values 
$$
 u_{\b,i-\frac12}^{\ast}= \frac 1{l_{\b,i+\frac12}}\sum_{\a= \a^-_{i-\frac12}(\b)  }^{\a^+_{i-\frac12}(\b)}  l_{\a,i-\frac12} u_{\a,i-\frac12}^{n}
$$
 will be used to compute the approximation of the velocity gradient   at $i+\frac12$.
  Clearly,  this may result in a local loss of accuracy, but the numerical results reported show that this has limited
 impact on the overall accuracy of the proposed method. 

\bigskip
 The discretization of the mass transfer term is defined as
\begin{eqnarray}
\label{eq:trasf_disc}
\mathcal{G}_{\a+\frac12,i+\frac12}^{n} &=&  \dsum_{\b=1}^{\a}   l_{\b,i+\frac12} \left(\left(h u_{\b} \,-\, \dsum_{\gamma=1}^N l_{\gamma}h u_{\gamma}\right)_{i+1}^n - \left(h u_{\b} \,-\, \dsum_{\gamma=1}^N l_{\gamma}h u_{\gamma}\right)_{i}^n\right),\nonumber
\end{eqnarray}
 where $ \left(h u_{\b} \,-\, \dsum_{\gamma=1}^N l_{\gamma}h u_{\gamma}\right)_{i}$ is the upwind value depending on the averaged velocity $u_{\b,i} = (u_{\b,i-\frac12}+u_{\b,i+\frac12})/2$. For the tracer equation, this term appears at the center of the control volume,  so that we set instead

\begin{eqnarray}
\label{trasf_disc}
G^n_{\a+\frac12,i} &=&  \dfrac{1}{\dxi}\dsum_{\b=1}^{\a} \left(  l_{\b,i+\frac12}h^n_{i+\frac12}u^n_{\b,i+\frac12} \,-\,l_{\b,i-\frac12}h^n_{i-\frac12}u^n_{\b,i-\frac12}   \right. \nonumber \\
&-& \left. l_{\b,i}\dsum_{\gamma=1}^{N_{i}}\left(  l_{\gamma,i+\frac12}h^n_{i+\frac12}u^n_{\gamma,i+\frac12} \,-\,l_{\gamma,i-\frac12}h^n_{i-\frac12}u^n_{\gamma,i-\frac12}    \right)\right).
\end{eqnarray}
 The above formulas are to be modified appropriately for cells in which
 $N_{i-\frac12} \neq N_{i+\frac12}, $ by summing all the contributions
 on the cell boundary with more layers that correspond to a given term 
 $ l_{\b,i\pm\frac12}h^n_{i\pm\frac12}u^n_{\b,i\pm\frac12} $ on the cell boundary with less layers, according to the definitions in the previous section.

 \begin{remark} The time discretization of the mass transfer terms could be easily turned into an implicit one, by taking instead
  \begin{eqnarray}
\label{eq:vel_discrete_theta_fullimp}
&&u_{\a,i+\frac12}^{n+1} + g\theta\frac{\dt} {\dxip}(\eta_{i+1}^{n+1} - \eta_{i}^{n+1}) \\
&&-\frac{\dt \theta}
{l_{\a,i+\frac12}\,h_{i+\frac12}^{n}}
\left(\gamma_{\a+\frac12,i+\frac12}^{n}\,\left(u_{\a+1,i+\frac12}^{n+1}-u_{\a,i+\frac12}^{n+1}\right) - \delta_{\a-\frac12,i+\frac12}^{n}\,\left(u_{\a,i+\frac12}^{n+1}-u_{\a-1,i+\frac12}^{n+1}\right) \right) \nonumber\\
&&= u_{\a,i+\frac12}^{n}
 + \dt  {\cal A}_{\a,i+\frac 12}^{u,n}   -   g(1-\theta)\frac{\dt}{\dxip}\left(\eta_{i+1}^{n} - \eta_{i}^{n}\right) \nonumber
 \\ 
&&+\frac{\dt (1-\theta)}{l_{\a,i+\frac12}\,h_{i+\frac12}^{n}}\left(\gamma_{\a+\frac12,i+\frac12}^{n}\,\left(u_{\a+1,i+\frac12}^{n}-u_{\a,i+\frac12}^{n}\right) 
- \delta_{\a-\frac12,i+\frac12}^{n}\,\left(u_{\a,i+\frac12}^{n}-u_{\a-1,i+\frac12}^{n}\right) \right), \nonumber
\end{eqnarray}
 where now 
 $$\gamma_{\a+\frac12,i+\frac12}^{n}=\dfrac{\nu_{\a+\frac12,i+\frac12}^{n}}{l_{\a+\frac12,i+\frac12}h_{i+\frac12}^{n}}
  + \frac{\mathcal{G}_{\a+\frac12,i+\frac12}^{n} }{2\dxip}
  \ \ \ \delta_{\a-\frac12,i+\frac12}^{n}=\dfrac{\nu_{\a-\frac12,i+\frac12}^{n}}{l_{\a-\frac12,i+\frac12}h_{i+\frac12}^{n}}
  - \frac{\mathcal{G}_{\a-\frac12,i+\frac12}^{n} }{2\dxip}.$$
  
\noindent This approach might be helpful to relax stability restrictions if large
  values of $\mathcal{G}_{\a+\frac12,i+\frac12}^{n} $ arise. In the implementation employed to
  obtain the numerical results of section \ref{tests}, however, only the discretization
  \eqref{eq:vel_discrete_theta} was applied so far.
  \end{remark}

 \bigskip
 At the bottom ($\a=1$) and at the free surface ($\a=N_{i+\frac12}$) layers,
 the viscous terms are modified by the friction and drag terms at
  $\Gamma_{1/2}$ and $\Gamma_{N_{i+\frac12}+1/2},$
 respectively. We have then

  \begin{eqnarray}
\label{eq:vel_discrete_theta_bottom}
&& u_{1,i+\frac12}^{n+1} +   g\theta\frac{\dt} {\dxip}(\eta_{i+1}^{n+1} - \eta_{i}^{n+1})
\\ 
&&-\frac{\dt \theta}{ l_{1,i+\frac12}\,h_{i+\frac12}^{n}}\left(\nu_{1+\frac12,i+\frac12}^{n}\,\dfrac{u_{2,i+\frac12}^{n+1}-u_{1,i+\frac12}^{n+1}}{l_{1+\frac12,i+\frac12}h_{i+\frac12}^{n}} 
-C^n_{f,i+\frac12}\abs{u_{1,i+\frac12}^n}u_{1,i+\frac12}^{n+1} \right) \nonumber      \\
 &&=u_{1,i+\frac12}^{n}
 + \dt  {\cal A}_{1,i+\frac 12}^{u,n} 
+\frac{\dt} {\dxip l_{1,i+\frac12}\,h_{i+\frac12}^{n}}\Delta \wtd{u}_{\frac32,i+\frac12}^{\,n}
  \mathcal{G}_{\frac32,i+\frac12}^{n} -g(1-\theta)\frac{\dt}{\dxip}\left(\eta_{i+1}^{n} - \eta_{i}^{n}\right)\nonumber \\
 &&+\frac{\dt(1-\theta) }{ l_{1,i+\frac12}\,h_{i+\frac12}^{n}} 
\left(\nu_{1+\frac12,i+\frac12}^{n}\,\dfrac{u_{2,i+\frac12}^{n}-u_{1,i+\frac12}^{n}}{l_{1+\frac12,i+\frac12}h_{i+\frac12}^{n}} 
-C^n_{f,i+\frac12}\abs{u_{1,i+\frac12}^n}u_{1,i+\frac12}^{n} \right), \nonumber
\end{eqnarray}

 \begin{eqnarray}
\label{eq:vel_discrete_theta_top}
&&u_{N_{i+\frac12},i+\frac12}^{n+1}   +  g\theta\frac{\dt} {\dxip}(\eta_{i+1}^{n+1} - \eta_{i}^{n+1}) \\
&& +\frac{\dt \theta}{l_{N_{i+\frac12},i+\frac12}\,h_{i+\frac12}^{n}}
\left(  \nu_{N_{i+\frac12}-\frac12,i+\frac12}^{n}\,\dfrac{u_{N_{i+\frac12},i+\frac12}^{n+1}-u_{N_{i+\frac12}-1,i+\frac12}^{n+1}}{l_{N_{i+\frac12}-\frac12,i+\frac12}h_{i+\frac12}^{n}}  + \wtd{C_w}_{,i+\frac12}^n u_{N_{i+\frac12}}^{n+1}\right) \nonumber\\
&&=  u_{N_{i+\frac12},i+\frac12}^{n}
 + \dt {\cal A}_{N_{i+\frac12},i+\frac 12}^{u,n} 
 +\frac{\dt \wtd{C_w}_{,i+\frac12}^n }{l_{N_{i+\frac12},i+\frac12}\,h_{i+\frac12}^{n}}
  \left(\theta u^{n+1}_{w,i+\frac12}
 +(1-\theta) \left(u^{n}_{w,i+\frac12}-u^{n}_{N_{i+\frac12},i+\frac12}\right)\right)
 \nonumber \\ 
&&-\frac{(1-\theta)\dt}{l_{N_{i+\frac12},i+\frac12}\,h_{i+\frac12}^{n}}  
\nu_{N_{i+\frac12}-\frac12,i+\frac12}^{n}\,\dfrac{u_{N_{i+\frac12},i+\frac12}^{n}-u_{N_{i+\frac12}-1,i+\frac12}^{n}}{l_{N_{i+\frac12}-\frac12,i+\frac12}h_{i+\frac12}^{n}} \nonumber\\
&&+\frac{\dt} {\dxip l_{N_{i+\frac12},i+\frac12}\,h_{i+\frac12}^{n}} \Delta \wtd{u}_{N_{i+\frac12}-\frac12,i+\frac12}^{\,n} \mathcal{G}_{N_{i+\frac12}-\frac12,i+\frac12}^{n} - g\left(1-\theta\right)\frac{\dt} {\dxip}(\eta_{i+1}^{n} - \eta_{i}^{n}).\nonumber
\end{eqnarray}
Notice that, 
in previous equation,  we define $ \wtd{C_w}^n = C_w|u_{w,i+\frac12}^n - u_{N_{i+\frac12},i+\frac12}^n|.$ 

\newpage

Replacing the  expressions for the velocities at time step $n+1$  into the continuity equation
yields a linear system whose unknowns are the values of the free surface $\eta_i^{n+1}$.  
This can be done rewriting  the discrete momentum equations in matrix notation  as in  \cite{casulli:1994},
after rescaling both sides of the equations by $l_{\a,i+\frac12}h^n_{i+\frac12}.$
We denote by $\bF_{i+\frac12}^n$ collects all the explicit terms and by $\bA^n_{i+\frac12} $
the matrix resulting from the discretization of the vertical diffusion terms.
  Since it is a tridiagonal, positive definite, diagonally dominant matrix, it is invertible and
  its inverse also has the same properties.
We also define
$$
\bU_{i+\frac12} = \left(\begin{array}{c}
u_{1,i+\frac12}\\
:\\
u_{\a,i+\frac12}\\
:\\
u_{N_{i+\frac12},i+\frac12},
\end{array}\right);\quad\quad \bH_{i+\frac12} = 
 \left(\begin{array}{c}
l_{1,i+\frac12}h_{i+\frac12}\\
:\\
l_{\a,i+\frac12}h_{i+\frac12}\\
:\\
l_{N_{i+\frac12},i+\frac12}h_{i+\frac12},
\end{array}\right).
$$
As a result, one can write
$$
\bA_{i+\frac12}^n \,\bU_{i+\frac12}^{n+1}\, = \,\bF_{i+\frac12}^n\, 
  - g\,\theta\dfrac{\dt}{\dxip}\left(\eta_{i+1}^{n+1}-\eta_{i}^{n+1}\right)\, \bH_{i+\frac12}^n
$$

  \begin{eqnarray}
  \label{eq:mom_eq_matrix}
 \bU_{i+\frac12}^{n+1}\, &=& \left ( \bA_{i+\frac12}^n\right)^{-1}\bF_{i+\frac12}^n\, \nonumber \\
  &-& g\,\theta\dfrac{\dt}{\dxip}\left(\eta_{i+1}^{n+1}-\eta_{i}^{n+1}\right)\,
  \left( \bA_{i+\frac12}^n \right )^{-1}\bH_{i+\frac12}^n.
\end{eqnarray}
 The discrete continuity equation is rewritten in this matrix notation as
$$
\begin{array}{rl}
\dxi\eta_{i}^{n+1} =&  \dxi\eta_{i}^{n} - \theta\dt\left(  (\bH_{i+\frac12}^n)^{T}\,\bU_{i+\frac12}^{n+1} - (\bH_{i-\frac12}^n)^{T}\,\bU_{i-\frac12}^{n+1} \right)\\[4mm]
-& (1-\theta)\dt\left(  (\bH_{i+\frac12}^n)^{T}\,\bU_{i+\frac12}^{n} - (\bH_{i-\frac12}^n)^{T}\,\bU_{i-\frac12}^{n} \right).
\end{array}
$$
  Substituting formally equation \eqref{eq:mom_eq_matrix}
in the continuity equation yields the tridiagonal system

\begin{eqnarray}
\dxi\eta_i^{n+1} &-& g\,\theta^2\dt^2\left(\left[\bH^{T}\bA^{-1}\bH\right]^n_{i+\frac12}
 \dfrac{\eta_{i+1}^{n+1} - \eta_i^{n+1}}{\dxip}  \right. \\
 &-& \left. \left[\bH^{T}\bA^{-1}\bH\right]^n_{i-\frac12} 
 \dfrac{\eta_{i}^{n+1} - \eta_{i-1}^{n+1}}{\dxim}\right) \nonumber \\
&=&  \dxi\eta_i^n - \theta\dt \left(\left[\bH^T \bA^{-1} \bG\right]^n_{i+\frac12} - \left[\bH^T \bA^{-1} \bG\right]^n_{i-\frac12} \right)\, \nonumber \\
&-&  (1-\theta) \dt \left((\bH_{i+\frac12}^n)^{T}\,\bU_{i+\frac12}^{n} - (\bH_{i-\frac12}^n)^{T}\,\bU_{i-\frac12}^{n}\right).
\nonumber
\end{eqnarray}
The new values of the free surface $\eta_i^{n+1}$ are computed by solving this system. The updated values $\eta^{n+1}_i$ are  then replaced in \eqref{eq:mom_eq_matrix} to obtain $u^{n+1}_{\a,i+\frac12}$.  

Finally, the  evolution equation for  $\r_\a$ is discretized as
\begin{eqnarray}\label{eq:evol_rho}
\dxi l_{\a,i}h_{i}^{n+1}\r_{\a,i}^{\,n+1} &=&  \dxi l_{\a,i}h_{i}^{n}\r_{\a,i}^{\,n} \\ 
&-& \dt\left( l_{\a,i+\frac12} h^n_{i+\frac12} \r_{\a,i+\frac12}^n u^{n+\theta}_{\a,i+\frac12}
 \,-\,l_{\a,i-\frac12} h_{i-\frac12}^n \r_{\a,i-\frac12}^nu^{n+\theta}_{\a,i-\frac12} \right) \nonumber \\
&+&  \dxi\dt\left(\r_{\a+\frac12,i}^{\,n} G_{\a+\frac12,i}^n - \r_{\a-\frac12,i}^{\,n} G_{\a-\frac12,i}^n\right)\, ,
\nonumber
\end{eqnarray}
where 
$   u_\a^{n+\theta} = \theta u_\a^{n+1} + (1-\theta)u_\a^n.$
The values $\r_{\a,i\pm\frac12}^n, \r_{\a\pm\frac12,i}^n$ can be defined by appropriate numerical fluxes.
Also the discretization of the tracer equation could be easily turned into an implicit one in the vertical
if required for stability reasons, by setting 
\begin{eqnarray}\label{eq:evol_rho_fimp}
\dxi l_{\a,i}h_{i}^{n+1}\r_{\a,i}^{\,n+1} 
&-& \theta\dxi\dt\left(\r_{\a+\frac12,i}^{\,n+1} G_{\a+\frac12,i}^n - \r_{\a-\frac12,i}^{\,n+1} G_{\a-\frac12,i}^n\right)
=  \dxi l_{\a,i}h_{i}^{n}\r_{\a,i}^{\,n} \nonumber\\ 
&-& \dt\left( l_{\a,i+\frac12} h^n_{i+\frac12} \r_{\a,i+\frac12}^n u^{n+\theta}_{\a,i+\frac12}
 \,-\,l_{\a,i-\frac12} h_{i-\frac12}^n \r_{\a,i-\frac12}^nu^{n+\theta}_{\a,i-\frac12} \right) \nonumber \\
&+&  (1-\theta)\dxi\dt\left(\r_{\a+\frac12,i}^{\,n} G_{\a+\frac12,i}^n - \r_{\a-\frac12,i}^{\,n} G_{\a-\frac12,i}^n\right).
\end{eqnarray}
Notice that, as in formula \eqref{trasf_disc}, the previous definitions
are to be modified appropriately for cells in which
 $N_{i-\frac12} \neq N_{i+\frac12}, $ by summing all the contributions
 on the cell boundary with more layers that correspond to a given term 
 $l_{\a,i\pm\frac12} h^n_{i+\pm\frac12} \r_{\a,i\pm\frac12}^n u^{n+\theta}_{\a,i\pm\frac12} $
 on the cell boundary with less layers,
 according to the definitions in the previous sections.

  It is also important to remark that, if $\r_{\a,i}^{\,n+1}=\r_{\a,i}^{\,n}=1 $
  is assumed in either   \eqref{eq:evol_rho},  \eqref{eq:evol_rho_fimp}, 
  as long as a consistent flux is employed for the definition of $\r_{\a,i\pm\frac12}^n, \r_{\a\pm\frac12,i}^n,$ 
  discretizations of the first
  equation in \eqref{eq:FinalModel} are obtained,
   which then summed over the whole set of layers
  $\a=1,...,N_{i} $ yield exactly formula \eqref{eq:mass_discrete_theta}. This implies
  that complete consistency with the discretization of the continuity equation is guaranteed.
  The importance of this property for the numerical approximation of free surface problems
   has been discussed extensively in \cite{gross:2002}.
 
 \subsection{Second order IMEX-ARK2 method}   

 A more accurate   time discretization can be achieved employing an IMplicit EXplicit  (IMEX) Additive Runge Kutta method (ARK)  \cite{kennedy:2003}.
 These techniques address the discretization of ODE systems that can be written as
 $\mathbf{y}' = \mathbf{f}{s}(\mathbf{y},t) +\mathbf{f}_{ns}(\mathbf{y},t), $ where the $s$ and $ns$ subscripts denote the stiff
 and non stiff components of the system, respectively. 
 In the case of system \eqref{eq:sistema3D_u}, the non stiff term would include the momentum advection and mass exchange terms,
 while the stiff term would include free surface gradients and vertical viscosity terms.
 A generic   $s-$stage IMEX-ARK method can be defined as follows. If $l_{\max}$ is the number of intermediate states of the method, then for $l=1,\dots, l_{\max}$: 
\begin{eqnarray}
\label{imex-ark}
\mathbf{u}^{(l)}=\mathbf{u}^n&+&\Delta t  \sum_{m=1}^{l-1} \bigg( a_{lm}\mathbf{f}_{ns}(\mathbf{u}^{(m)},t+c_m \Delta t)+\tilde{a}_{lm}\mathbf{f}_{s}(\mathbf{u}^{(m)},t+c_m\Delta t) \bigg) + \nonumber\\ &+&\Delta t \, \tilde{a}_{ll} \, \mathbf{f}_{s}(\mathbf{u}^{(l)},t+c_i\Delta t), 
\end{eqnarray}
Finally, $u^{n+1}$ is computed:
$$
\mathbf{u}^{n+1}=\mathbf{u}^n+\Delta t \sum_{l=1}^{l_{\max}}b_{l}(\mathbf{f}_{ns}(\mathbf{u}^{(l)},t+c_l \Delta t)+\mathbf{f}_{s}(\mathbf{u}^{(l)},t+c_l\Delta t)).
$$
Coefficients $a_{lm}, \tilde{a}_{lm}, c_l$ and $b_l$  are determined so that the method is consistent of a given order. In addition to the order conditions specific to each sub-method, the coefficients should respect coupling conditions.
 Here, we will use  a specific second order method, whose coefficients are presented in  the Butcher tableaux. See tables \ref{ark2_butch_e} and \ref{ark2_butch_i} for the explicit and implicit method, respectively.
 The coefficients of  the explicit method were proposed in
\cite{giraldo:2013}, while the implicit  method, also employed in the same paper,
  coincides indeed with the TR-BDF2 method proposed in
\cite{bank:1985}, \cite{hosea:1996} and applied to the shallow water 
and Euler equations in \cite{tumolo:2015}.

\begin{table}[!h]
\begin{center}
\begin{tabular}{c|ccc}
0 & 0 & &  \\
$2 \mp \sqrt{2}$ & $2 \mp \sqrt{2}$ & 0 & \\
1 & $1-(3+2\sqrt{2})/6$ & $ (3+2\sqrt{2})/6$ & 0 \\
\hline
 & $\pm \frac{1}{2\sqrt{2}}$ & $\pm \frac{1}{2\sqrt{2}}$ & $1 \mp \frac{1}{\sqrt{2}}$
 \end{tabular}
 \end{center}
\caption{Butcher tableaux of the explicit ARK2 method}
 \label{ark2_butch_e}
 \end{table}

\begin{table}[h!]
\begin{center}
\begin{tabular}{c|ccc}
0 & 0 & & \\
$2 \mp \sqrt{2}$ &  $1 \mp \frac{1}{\sqrt{2}}$ &  $1 \mp \frac{1}{\sqrt{2}}$ & \\
1 & $\pm \frac{1}{2\sqrt{2}}$ & $\pm \frac{1}{2\sqrt{2}}$ &1 $\mp \frac{1}{2\sqrt{2}}$ \\
\hline
 & $\pm \frac{1}{2\sqrt{2}}$ & $\pm \frac{1}{2\sqrt{2}}$ & $1 \mp \frac{1}{\sqrt{2}}$
 \end{tabular}
 \end{center}
\caption{Butcher tableaux of the implicit-ARK2 method}
 \label{ark2_butch_i}
 \end{table}

While a straightforward application of \eqref{imex-ark} is certainly possible, we will outline here a more
efficient way to implement this method to the discretization
of equations \eqref{eq:sistema3D_u}, that mimics what done above for the   simpler $\theta-$method.
For the first stage, we define $ \eta_{i}^{n,1}=\eta_{i}^{n}, $   $u_{\a,i+\frac12}^{n,1}=u_{\a,i+\frac12}^{n}, $ 
and $\r_{\a,i}^{n,1}=\r_{\a,i}^{n}, $ 
respectively.
For the second stage, we get for the continuity equation
 \begin{eqnarray}
 \label{eq:mass_discrete_ark2}
&&\dxi \eta_{i}^{n,2}+ \,\tilde a_{22}\dt 
\left( \dsum_{\b=1}^{N_{i+\frac12}}
 l_{\b,i+\frac12} h^n_{i+\frac12} u_{\b,i+\frac12}^{n,2} \, - \, 
 \dsum_{\b=1}^{N_{i-\frac12}} l_{\b,i-\frac12} h^n_{i-\frac12}u_{\b,i-\frac12}^{n,2}\right) \nonumber \\
&&=   \dxi\eta_{i}^{n} -\,\tilde a_{21}\dt 
\left( \dsum_{\b=1}^{N_{i+\frac12}}l_{\b,i+\frac12} h^n_{i+\frac12} u_{\b,i+\frac12}^{n,1} \,
 - \,   \dsum_{\b=1}^{N_{i-\frac12}}l_{\b,i-\frac12} h^n_{i-\frac12}u_{\b,i-\frac12}^{n,1}\right) \nonumber
\end{eqnarray} 
and for the momentum equations 
  
 \begin{eqnarray}
\label{eq:vel_discrete_ark2}
&&u_{\a,i+\frac12}^{n,2} + g\tilde a_{22}\frac{\dt} {\dxip}(\eta_{i+1}^{n,2} - \eta_{i}^{n,2})\nonumber \\ 
&&-\frac{\dt \tilde a_{22}} {l_{\a,i+\frac12}\,h_{i+\frac12}^{n}}  \left(\nu_{\a+\frac12,i+\frac12}^{n}\,\dfrac{u_{\a+1,i+\frac12}^{n,2}
-u_{\a,i+\frac12}^{n,2}}{l_{\a+\frac12}h_{i+\frac12}^{n}} - \nu_{\a-\frac12,i+\frac12}^{n}\,\dfrac{u_{\a,i+\frac12}^{n,2}
-u_{\a-1,i+\frac12}^{n,2}}{l_{\a-\frac12}h_{i+\frac12}^{n}} \right) \nonumber\\[2mm]
&&=u_{\a,i+\frac12}^{n}
  + \dt a_{21}{\cal F}_{\a,i+\frac 12}^{\, n,1}   + \dt \tilde a_{21} {\cal I}^{\, n,1}_{\a,i+\frac12}
\end{eqnarray}
for $\a=1,...,N,$ with the appropriate corrections for the top and bottom layers, respectively.
Here we define
$${\cal F}_{\a,i+\frac 12}^{\, n,j}= {\cal A}_{\a,i+\frac 12}^{u,n,j}+\dfrac{1}{\dxip l_{\a,i+\frac12}\,h_{i+\frac12}^{n}}\left(\Delta \wtd{u}_{\a+\frac12,i+\frac12}^{\,n,j} \mathcal{G}_{\a+\frac12,i+\frac12}^{n,j} 
+\Delta \wtd{u}_{\a-\frac12,i+\frac12}^{\,n,j} \mathcal{G}_{\a-\frac12,i+\frac12}^{n,j}\right) $$
and
$$\begin{array}{l}
{\cal I}^{\, n,j}_{\a,i+\frac12} = -  \dfrac{g}{\dxip}\left(\eta_{i+1}^{n,j} - \eta_{i}^{n,j}\right)\\
+\dfrac{1} {l_{\a,i+\frac12}\,h_{i+\frac12}^{n}}  
\left(\nu_{\a+\frac12,i+\frac12}^{n}\,\dfrac{u_{\a+1,i+\frac12}^{n,j}-u_{\a,i+\frac12}^{n,j}}{l_{\a+\frac12}h_{i+\frac12}^{n}} 
- \nu_{\a-\frac12,i+\frac12}^{n}\,\dfrac{u_{\a,i+\frac12}^{n,j}-u_{\a-1,i+\frac12}^{n,j}}{l_{\a-\frac12}h_{i+\frac12}^{n}} \right),\nonumber
\end{array}$$
and all the other symbols have the same interpretation as in the presentation of the $\theta-$ method.
It can be noticed that, again, the dependency on $h$ has been frozen at time level $n$ in order to
 avoid solving a nonlinear system at each timestep. As shown in \cite{boittin:2015},
 \cite{casulli:1994}, \cite{tumolo:2015},  this does not degrade the accuracy of the method. Also the dependency on time of the vertical viscosity is frozen at time level $n.$
 The same will be done for both kinds of coefficients also in the third stage of the method.
 As in the previous discussion, the above discrete equations can be
 rewritten in vector notation as
 
   \begin{equation}
  \label{eq:mom_eq_matrix_ark2_2}
 \bU_{i+\frac12}^{n,2}\,=  \left ( \bA_{i+\frac12}^n\right)^{-1}\bF_{i+\frac12}^{n,1}\,
  -  g\,\tilde a_{22}\dfrac{\dt}{\dxip}\left(\eta_{i+1}^{n,2}-\eta_{i}^{n,2}\right)\,\left( \bA_{i+\frac12}^n \right )^{-1}\bH_{i+\frac12}^n,
\end{equation}
  where now  $\bF_{i+\frac12}^1$ has components given by 
 $$ l_{\a,i+\frac12}\,h_{i+\frac12}^{n}\left(u_{\a,i+\frac12}^{n}+ \dt a_{21}{\cal F}_{\a,i+\frac 12}^{\, n,1}\right. + \left.\dt \tilde{a}_{21}{\cal I}_{\a,i+\frac 12}^{\, n,1}\right).$$
The discrete continuity equation is rewritten in this matrix notation as
$$
\begin{array}{rl}
\dxi\eta_{i}^{n,2} =&  \dxi\eta_{i}^{n,2} 
- \tilde a_{22}\dt\left(  (\bH_{i+\frac12}^n)^{T}\,\bU_{i+\frac12}^{n,2} - (\bH_{i-\frac12}^n)^{T}\,\bU_{i-\frac12}^{n,2} \right)\\[4mm]
-& \tilde a_{21}\dt\left(  (\bH_{i+\frac12}^n)^{T}\,\bU_{i+\frac12}^{n} - (\bH_{i-\frac12}^n)^{T}\,\bU_{i-\frac12}^{n} \right).
\end{array}
$$
  Substituting formally equation \eqref{eq:mom_eq_matrix_ark2_2}
in the momentum equation yields the tridiagonal system
$$
\begin{array}{rl}
\dxi\eta_i^{n,2} - & g\, \tilde a_{22}^2\dt^2\left(\left[\bH^{T}\bA^{-1}\bH\right]^n_{i+\frac12}
 \dfrac{\eta_{i+1}^{n,2} - \eta_i^{n,2}}{\dxip} - \left[\bH^{T}\bA^{-1}\bH\right]^n_{i-\frac12} 
 \dfrac{\eta_{i}^{n,2} - \eta_{i-1}^{n,2}}{\dxim}\right) \\[4mm]
= & \dxi\eta_i^n -  \tilde a_{22}\dt\left(\left[\bH^T \bA^{-1} \bF^1\right]^n_{i+\frac12} - \left[\bH^T \bA^{-1} \bF^1\right]^n_{i-\frac12} \right)\\[4mm]
- & \tilde a_{21}\dt\left((\bH_{i+\frac12}^n)^{T}\,\bU_{i+\frac12}^{n} - (\bH_{i-\frac12}^n)^{T}\,\bU_{i-\frac12}^{n}\right).
\end{array}$$
The new values of the free surface $\eta_i^{n,2}$ are computed by solving this system  and
they are replaced in \eqref{eq:mom_eq_matrix_ark2_2} to
 find $u^{n,2}_{\a,i+\frac12}. $ The  evolution equation for  $\r_\a$ is then discretized as
\begin{eqnarray}\label{eq:evol_rho_imex1}
\dxi l_{\a,i}h_{i}^{n,2}\r_{\a,i}^{\,n,2} &=&  \dxi l_{\a,i}h_{i}^{n}\r_{\a,i}^{\,n,1}\\ 
&-& a_{21}\dt\left( l_{\a,i+\frac12} h^{n}_{i+\frac12} \r_{\a,i+\frac12}^{n,1} u^{*,2}_{\a,i+\frac12}
 \,-\,l_{\a,i-\frac12} h_{i-\frac12}^{n} \r_{\a,i-\frac12}^{n,1}u^{*,2}_{\a,i-\frac12} \right) \nonumber \\
&+&  a_{21}\dxi\dt\left(\r_{\a+\frac12,i}^{n,1} G_{\a+\frac12,i}^{n,1} - \r_{\a-\frac12,i}^{n,1} G_{\a-\frac12,i}^{n,1}\right)\, ,
\nonumber
\end{eqnarray}
where now
$   u_\a^{*,2} =  \tilde a_{22} u_\a^{n,2} + \tilde a_{21}u_\a^{n,1}.$

 The last stage of the IMEX-ARK2 method can then be written in vector notation as
  \begin{eqnarray}
  \label{eq:mom_eq_matrix_ark2_3}
  \bU_{i+\frac12}^{n,3}\, &=& \left ( \bA_{i+\frac12}^n\right)^{-1}\bF_{i+\frac12}^{\, n,2}\, \nonumber \\
  &-& g\,\tilde a_{33}\dfrac{\dt}{\dxip}\left(\eta_{i+1}^{n,3}-\eta_{i}^{n,3}\right)\,\left( \bA_{i+\frac12}^n \right )^{-1}\bH_{i+\frac12}^n,
\end{eqnarray}
  where now  $\bF_{i+\frac12}^{\, n,2}$ has components given by 
 $$ l_{\a,i+\frac12}\,h_{i+\frac12}^{n}\left(u_{\a,i+\frac12}^{n}+\dt a_{31}{\cal F}_{\a,i+\frac 12}^{\, n,1}\right. + \left.\dt a_{32}{\cal F}_{\a,i+\frac 12}^{\, n,2}+\dt \tilde{a}_{31}{\cal I}_{\a,i+\frac 12}^{\, n,1} + \dt \tilde{a}_{32}{\cal I}_{\a,i+\frac 12}^{\, n,2}\right).$$
The discrete continuity equation is rewritten in this matrix notation as
$$
\begin{array}{rl}
\dxi\eta_{i}^{n,3} =&  \dxi\eta_{i}^{n,3} 
- \tilde a_{33}\dt\left(  (\bH_{i+\frac12}^n)^{T}\,\bU_{i+\frac12}^{n,3} - (\bH_{i-\frac12}^n)^{t}\,\bU_{i-\frac12}^{n,3} \right)\\[4mm]
-& \tilde a_{31}\dt\left(  (\bH_{i+\frac12}^n)^{T}\,\bU_{i+\frac12}^{1} - (\bH_{i-\frac12}^n)^{T}\,\bU_{i-\frac12}^{1} \right) \\[4mm]
-& \tilde a_{32}\dt\left(  (\bH_{i+\frac12}^n)^{T}\,\bU_{i+\frac12}^{2} - (\bH_{i-\frac12}^n)^{T}\,\bU_{i-\frac12}^{2} \right). 
\end{array}
$$
As a result, substitution of \eqref{eq:mom_eq_matrix_ark2_3} into the third stage of the continuity equation yields
the tridiagonal system
$$
\dxi\eta_i^{n,3} -  g\, \tilde a_{33}^2\dt^2\left(\left[\bH^{T}\bA^{-1}\bH\right]^n_{i+\frac12}
 \dfrac{\eta_{i+1}^{n,3} - \eta_i^{n,3}}{\dxip} - \left[\bH^{T}\bA^{-1}\bH\right]^n_{i-\frac12} 
 \dfrac{\eta_{i}^{n,3} - \eta_{i-1}^{n,3}}{\dxim}\right) = {\cal E}_i,$$
where now all the explicit terms have been collected in ${\cal E}_i.$
The new values of the free surface $\eta_i^{n,3}$ are computed by solving this system  and
they are replaced in \eqref{eq:mom_eq_matrix_ark2_3} to
 find $u^{n,3}_{\a,i+\frac12}.$ 
 The tracer density is then updated as
 \begin{eqnarray}\label{eq:evol_rho_imex2}
\dxi l_{\a,i}h_{i}^{n,3}\r_{\a,i}^{\,n,3} &=&  \dxi l_{\a,i}h_{i}^{n,2}\r_{\a,i}^{\,n,2}\\ 
&-& a_{32}\dt\left( l_{\a,i+\frac12} h^{n}_{i+\frac12} \r_{\a,i+\frac12}^{n} u^{*,3}_{\a,i+\frac12}
 \,-\,l_{\a,i-\frac12} h_{i-\frac12}^{n} \r_{\a,i-\frac12}^{n,2}u^{*,3}_{\a,i-\frac12} \right) \nonumber \\
&+&  a_{32}\dxi\dt\left(\r_{\a+\frac12,i}^{n,2} G_{\a+\frac12,i}^{n,2} - \r_{\a-\frac12,i}^{n,1} G_{\a-\frac12,i}^{n,2}\right)\nonumber\\
&-& a_{31}\dt\left( l_{\a,i+\frac12} h^{n}_{i+\frac12} \r_{\a,i+\frac12}^{n,1} u^{*,2}_{\a,i+\frac12}
 \,-\,l_{\a,i-\frac12} h_{i-\frac12}^{n} \r_{\a,i-\frac12}^{n,1}u^{*,2}_{\a,i-\frac12} \right) \nonumber \\
&+&  a_{31}\dxi\dt\left(\r_{\a+\frac12,i}^{n,1} G_{\a+\frac12,i}^{n,1} - \r_{\a-\frac12,i}^{n,1} G_{\a-\frac12,i}^{n,1}\right)
\, ,
\nonumber
\end{eqnarray}
where now
$   u_\a^{*,3} =  \tilde a_{32} u_\a^{n,3} + \tilde a_{31}u_\a^{n,2}.$

 The final assembly of the solution at time level $n+1$
 has then the form
 \begin{eqnarray}\label{eq:mass_discrete_ark2_final}
\dxi \eta_{i}^{n+1}&=& \dxi \eta_{i}^{n}\\
&-& \,\dt  \dsum_{j=1}^{3}\tilde b_{j} \left(  \dsum_{\b=1}^{N_{i+\frac12}}
l_{\b,i+\frac12} h^n_{i+\frac12} u_{\b,i+\frac12}^{n,j} \, - \, \dsum_{\b=1}^{N_{i-\frac12}} l_{\b,i-\frac12} h^n_{i-\frac12}u_{\b,i-\frac12}^{n,j}\right)  
\nonumber 
\end{eqnarray}
for the continuity equation,  
 
 \begin{eqnarray}
\label{eq:vel_discrete_ark2_final}
&&u_{\a,i+\frac12}^{n+1} =u_{\a,i+\frac12}^{n} +  \dt \dsum_{j=1}^{3}\left(\tilde{b}_{j} {\cal I}_{\a,i+\frac 12}^{n,j}   + b_{j} {\cal F}_{\a,i+\frac 12}^{n,j} \right)
\end{eqnarray}
for the momentum equations for $\a=1,...,N_{i+\frac 12},$ with the appropriate 
  corrections for the top and bottom layers, respectively, and
  
  \begin{eqnarray}\label{eq:evol_rho_imex_final}
\dxi l_{\a,i}h_{i}^{n+1}\r_{\a,i}^{\,n+1} &=&  \dxi l_{\a,i}h_{i}^{n}\r_{\a,i}^{\,n} 
- \dt \dsum_{j=1}^{3} b_{j} \left( l_{\a,i+\frac12} h^{n}_{i+\frac12} \r_{\a,i+\frac12}^{n,j} u^{n,j}_{\a,i+\frac12} \right .
 \nonumber \\
 &-&\left. l_{\a,i-\frac12} h_{i-\frac12}^{n} \r_{\a,i-\frac12}^{n,j}u^{n,j}_{\a,i-\frac12} \right)\\
&+&   \dxi\dt  \dsum_{j=1}^{3} b_{j} \left(\r_{\a+\frac12,i}^{n,j} G_{\a+\frac12,i}^{n,j}
 - \r_{\a-\frac12,i}^{n,j} G_{\a-\frac12,i}^{n,j}\right). \nonumber
\end{eqnarray}

  Notice that, also in this case, consistency with the discrete continuity equation in the sense of \cite{gross:2002} is guaranteed and an implicit treatment of the vertical advection term would be feasible with the same procedure
  outlined above for the $\theta-$method.
  Furthermore, the two linear systems that must be solved for each time step have identical structure and matrices that only differ by a constant   factor, thanks to the freezing of their coefficients at time level $n.$ This implies that,
  recomputing their entries does not entail a major overhead. It was shown in \cite{tumolo:2015} that, while apparently more   costly than the simpler $\theta-$method, this procedure leads indeed to an increase in efficiency by significantly increasing the accuracy that can be achieved with a given time step.

 \bigskip
 \section{Numerical results}
 \label{tests}
 In this section, we describe the results of several numerical experiments that were performed in order to investigate the accuracy and efficiency of the proposed  methods. In particular,  the potential loss of accuracy when reducing the number of vertical layers is investigated in each test,  as well as the reduction of the number of degrees of freedom of the system achieved by simplifying the vertical discretization in certain areas of the domain.
 
 We define the maximum Courant number associated to the velocity $C_{vel}$ and to the celerity $C_{cel}$  as 
 \begin{eqnarray}\label{eq:cfl}
C_{vel} &=& \max_{1 \leq i \leq M} \max_{1 \leq \a \leq N} \abs{u_{\a,{i+\frac12}}}\,\dfrac{\dt}{\dxi}; \nonumber \\
 C_{cel} &=& \max_{1 \leq i \leq M}\,\max_{1 \leq \a \leq N} \left(\abs{u_{\a,{i+\frac12}}} + \sqrt{g\,h_i}\right)\,\dfrac{\dt}{\dxi}.
 \end{eqnarray}
 In order to evaluate the accuracy of the semi-implicit schemes, we compute the relative errors between the computed solution and a reference solution. We denote by $Err_\eta\,[\,l_2\,]$ and  $Err_\eta\,[l_\infty\,]$ the relative error for the free surface when considering the usual $l_2 $ and $l_\infty$ norm, respectively. For the velocity we define
 \begin{equation}\label{eq:error}
 \begin{array}{l}
 Err^2_u = \left(\dfrac{\sum_{\a = 1}^N \sum_{i = 1}^M \abs{u_{\a,{i+\frac12}}-u_{\a,{i+\frac12}}^{ref}}^2\dxi h_{\a,i}  }{\sum_{\a = 1}^N \sum_{i = 1}^M \abs{u_{\a,{i+\frac12}}^{ref}}^2\dxi h_{\a,i} }\right)^{1/2};\\
 \\
 Err_u^{\infty} = \dfrac{\max_{\a}\max_{i} \abs{u_{\a,{i+\frac12}}-u_{\a,{i+\frac12}}^{ref}} }{\max_{\a}\max_{i} \abs{u_{\a,{i+\frac12}}^{ref}}  },
 \end{array}
 \end{equation}
 where $u^{ref}$ denotes the reference solution.
We consider as a reference solution the one computed by using an explicit third order Runge Kutta method with a maximum value for the celerity Courant number of $0.1$. Therefore, for the explicit scheme the Courant number is fixed and we consider an adaptive time step.

   \subsection{Free oscillations in a closed basin}
   \label{free_osc}
We consider here a subcritical flow in a closed domain of length  $L = 10$ km. The bottom topography is given by the Gaussian function 
 $$b(x) = 4\,e^{-(x-x_0)^2/\sigma^2},$$
 where $x_0 = 5000$ m and $\sigma = 0.1\,L$.   At the initial time  the flow is at rest and we take as initial free surface  
profile $\eta_0(x) = 10   + a x, $ where $a $ is chosen so that
 the water height is $h=10$ m at $x=0$  and $h=11$ m at $x=10$ km. We simulate the resulting oscillations  until $t=10800$ s ($3$ h). All the simulations are performed by using $10$ layers in the multilayer code and a uniform space discretization step $\dx = 50$ m. The friction coefficient $C_f$ is defined by (\ref{eq_def_cf}) with $\D z_r = h_1$ ($h_1 = l_1 h$), $\D z_0 = 3.3 \times10^{-5}$ and $\kappa = 0.41$. The wind drag is defined by  the coefficient value $C_w = 1.2\times 10^{-6}$  and we set a constant wind velocity $\vu_{w}=-1$ $m/s$.

\begin{figure}
 \hspace{-1.4cm} \includegraphics[width=1.2\textwidth]{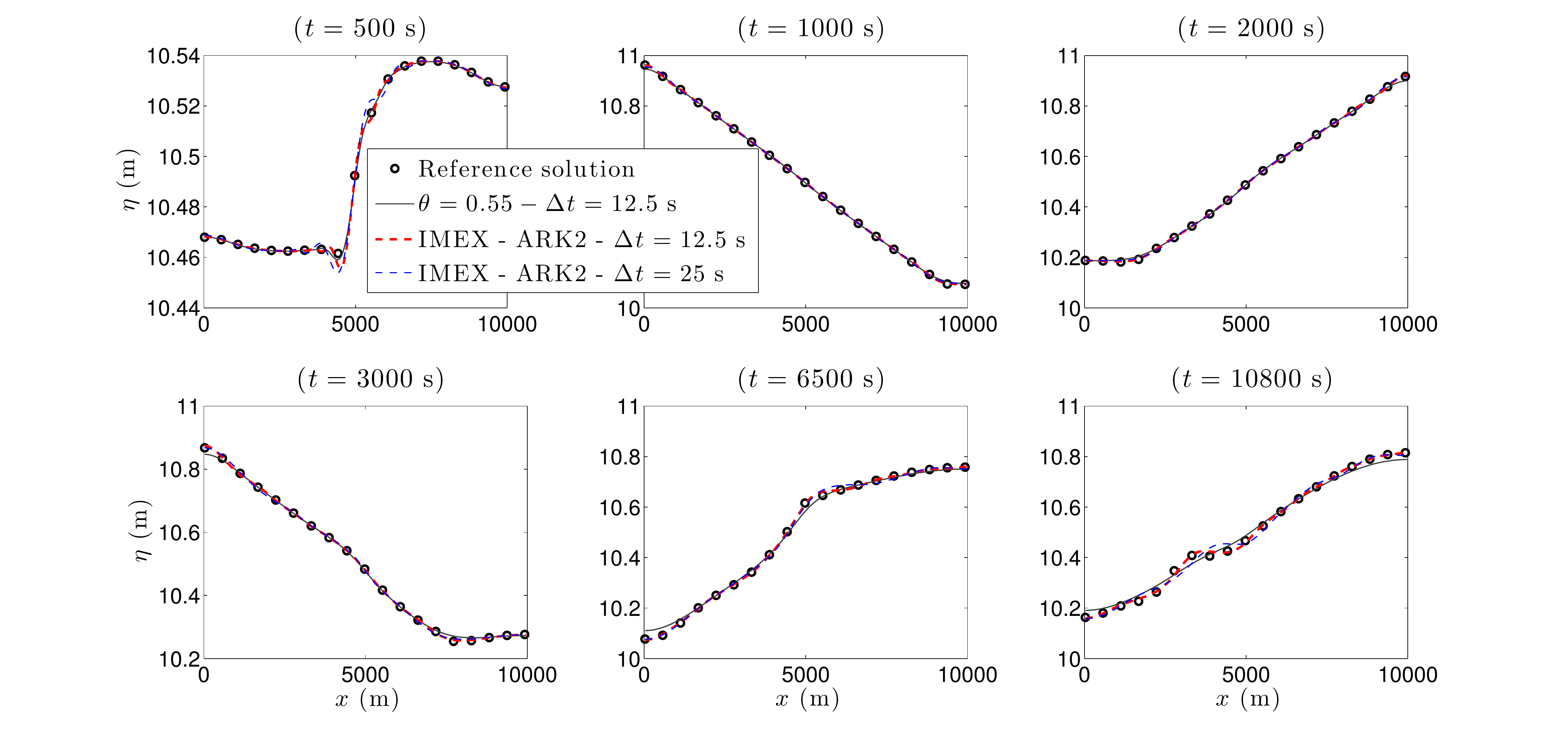}
  \caption{\label{fig:oscilation} \footnotesize \textit{Free surface profile at different times 
  by using the semi-implicit methods (color lines) and the reference solution (black circles) computed with the explicit scheme in the free oscillations test.}}
 \end{figure}

In figure \ref{fig:oscilation} we show  free surface profiles at different times until the final time, as computed  with the semi implicit methods described in section \ref{si_td}. The $\theta$-method and the IMEX-ARK2 are very close to the reference solution. However, the IMEX-ARK2 captures the shape of the free surface slightly better that the $\theta$-method when considering the same time step. By using the implicitness parameter $\theta = 0.55$ and the IMEX-ARK2 with $\dt = 12.5$ or $25$ s, we get a difference in the free surface of approximately 3 cm at the final time. In table \ref{tab:tabla1} we report the corresponding relative errors and the maximum Courant number achieved   by \eqref{eq:cfl}-\eqref{eq:error}, at   time $t = 10000$ s. We see that the IMEX-ARK2  method 
slightly improves the results  with respect to the $\theta$-method.

Even though it is hard to make a rigorous efficiency comparison in the framework of our preliminary implementation, for the subcritical regime the semi-implicit methods are much turn out to be more efficient than the explicit one. Actually, the computing time required to get the 3 hours of simulation (on a Mac Mini with Intel$\textsuperscript{\textregistered}$Core$\texttrademark$ i7-4578U and 16 GB of RAM) is approximately $12$ s for the explicit scheme using the Courant number $C_{cel} = 0.9$ ($103$ s for the reference solution), while it is approximately $1.64$ s (3.83 s) for the $\theta$-method (IMEX-ARK2) with $\dt = 12.5$ s. This time is $0.82$ s ($1.92$ s) with $\dt = 25$s and $0.4$ s ($0.97$ s) when considering the time step $\dt = 50$ s.  
 
\begin{figure}
 \hspace{-1.4cm} \includegraphics[width=1.2\textwidth]{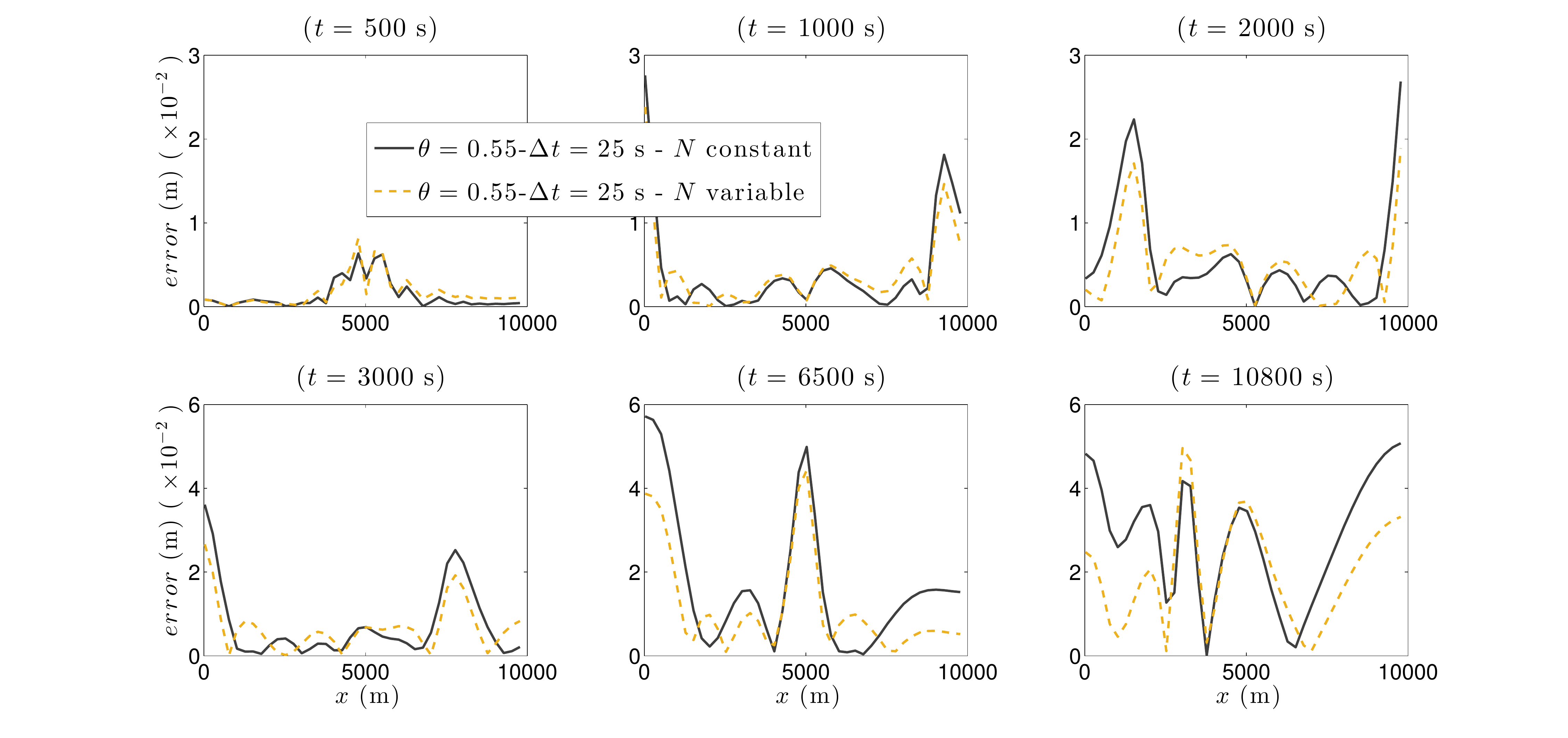}
  \caption{\label{fig:oscilation_nvar} \footnotesize \textit{Absolute errors for the free surface  at different times in the free oscillations test, obtained with the $\theta$-method ($\theta = 0.55$ and $\dt = 25$ s) and either   $10$ layers in the whole domain (solid black line) or a single layer in the first half of the domain  only   (dashed yellow line).}}
   \end{figure}

We then compare results obtained with a fixed and variable number of vertical layers.
Figure \ref{fig:oscilation_nvar} shows the absolute error for the free surface by using the $\theta$-method with $\theta=0.55$ and $\dt =25$ s,  as computed using either $N=10 $ layers throughout the domain or  considering
\begin{equation}\label{eq:nvar}
N = \left\{\begin{array}{ccl}
10 & &\text{if}\quad x\leq5000,\\
1 & &\text{otherwise}.
\end{array}\right.
\end{equation}

\begin{figure}
 \hspace{-1.4cm} \includegraphics[width=1.2\textwidth]{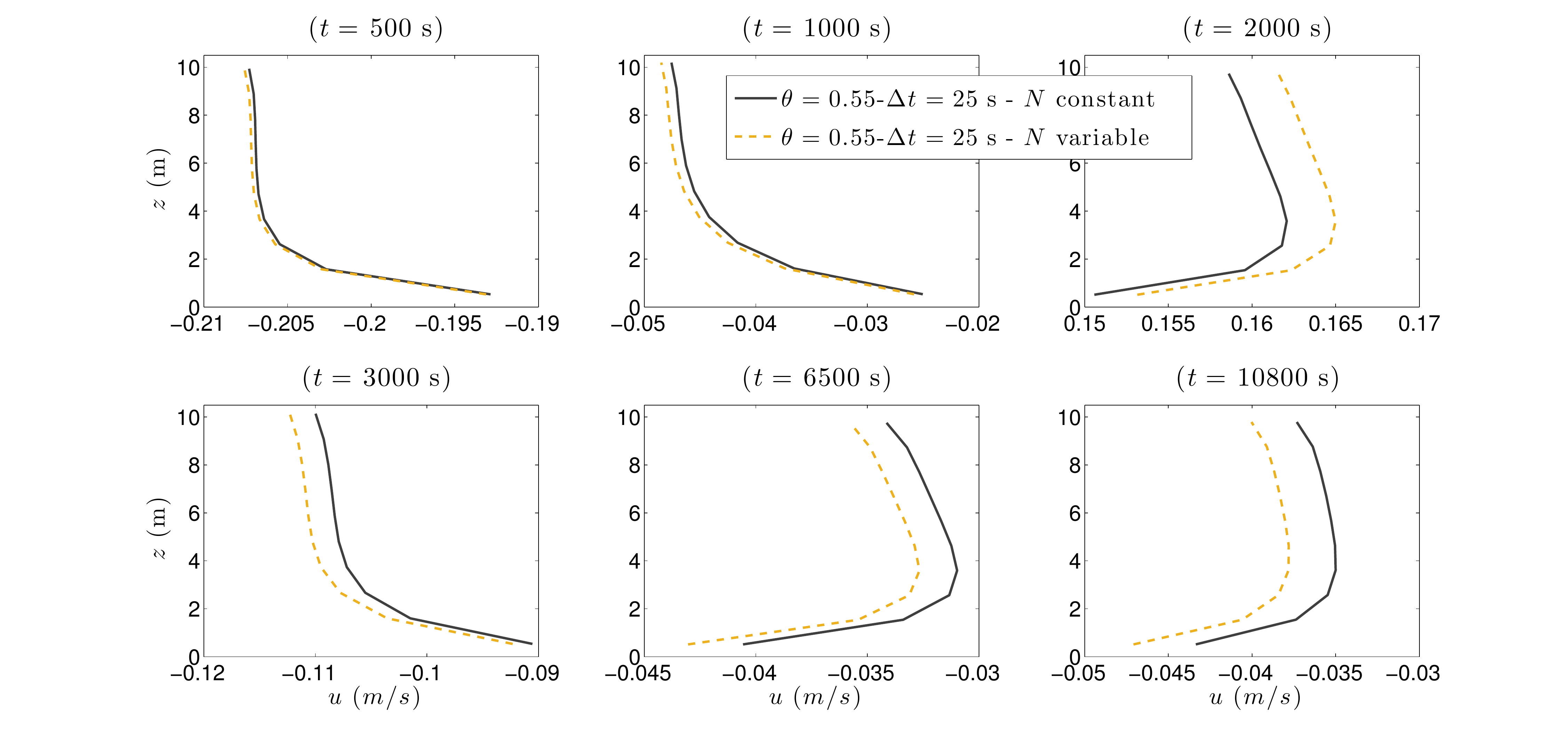}
  \caption{\label{fig:oscilation_vertical_nvar} \footnotesize \textit{Vertical profiles of horizontal velocity in the free oscillations test, obtained with the $\theta$-method  ($\theta = 0.55$ and $\dt = 25$ s) 
  and either   $10$ layers in the whole domain (solid black line) or a single layer in the first half of the domain  only   (dashed yellow line). Profiles  are taken at the point $x = 2475$ m   and times $t = 500,\,1000,\,2000,\,3000,\,6500,\,10800$ s.}}
 \end{figure}

Similar results are obtained if the time step is $\dt =12.5$ s. We see that usually the difference between the constant and variable layer cases computed by the semi-implicit method is of the order of 0.1\% of the solution values (absolute error 1 cm), while the number of degrees of freedom of the multilayer system is significantly reduced (from 2210 to 1310). Moreover, figure \ref{fig:oscilation_vertical_nvar} shows the vertical profiles of horizontal velocity at the point $x=2475$ m, as computed by the semi-implicit method with a constant and variable  number of layers (see \eqref{eq:nvar}).

\begin{table}
 \begin{center}
 \begin{tabular}{cccccccc}
 \hline\\
SI-method & $\dt$ (s) & $C_{vel}$  & $C_{cel}$ & Err$_{\eta} $ [$l_2/l_\infty$]  & Err$_{u}$ [$l_2/l_\infty$] \\ 
 & & &  &   ${\small(\times10^{-3})}$ &  ${\small(\times10^{-1})}$\\ \hline\\
$\theta=0.55$ & 12.5 & 0.18 & 2.62 & 1.6/3.2  & 0.9/1.5\\
IMEX-ARK2 & 12.5 & 0.18 & 2.62 & 0.6/2.0  & 0.4/0.6\\
$\theta=0.55$ & 25 & 0.34 & 5.24 & 2.6/5.4  & 1.3/1.7\\
IMEX-ARK2 & 25 & 0.34 & 5.24 & 0.9/2.2  & 1.2/1.7\\
$\theta=0.52$ & 50 & 0.7 & 10.48 & 3.1/6.3 & 1.6/1.5\\
$\theta=0.55$ & 50 & 0.68 & 10.47 & 3.9/7.7  & 2.2/2.0\\
IMEX-ARK2 & 50 & 0.69 & 10.48 & 2.4/5.2  & 1.4/1.7\\
\hline\end{tabular}
 \caption{\footnotesize \it{ Relative errors and Courant numbers achieved by using semi-implicit methods 
 in the free oscillations test at $t= 10000$ s.}}
\label{tab:tabla1}
  \end{center}
 \end{table}

    \subsection{Steady subcritical flow over a peak with friction}
   \label{subcritical}
In this test, a steady flow in the subcritical regime is considered, as done for example in \cite{rosatti:2011}.
The length of the domain is $L = 50$ m, and the bottom bathymetry is given by the function
  \begin{equation}
  b(x) = 0.05 - 0.001x + \left\{\begin{matrix}
  2 \cos^2 \left(\dfrac{\pi x}{10}\right), & \abs{x} < 5;\\[3mm]
  0 & \text{otherwise.}
  \end{matrix}\right. 
  \end{equation} 
 The initial conditions are given by $\eta_0(x) = 5$ m and $q_0(x) = 4.42$ $m^2\,s^{-1}$ and subcritical boundary conditions are considered. The same values of discharge and free surface are used for the upstream condition $q(-25,t)$, and the downstream one $\eta(L,t)$. We take a uniform space discretization step $\dx = 0.25$ m and the same values for the turbulent viscosity and bottom friction as in the previous test, while the wind stress is not taken into account in this case. 
   \begin{figure}[!h]\centering
\includegraphics[width=0.6\textwidth]{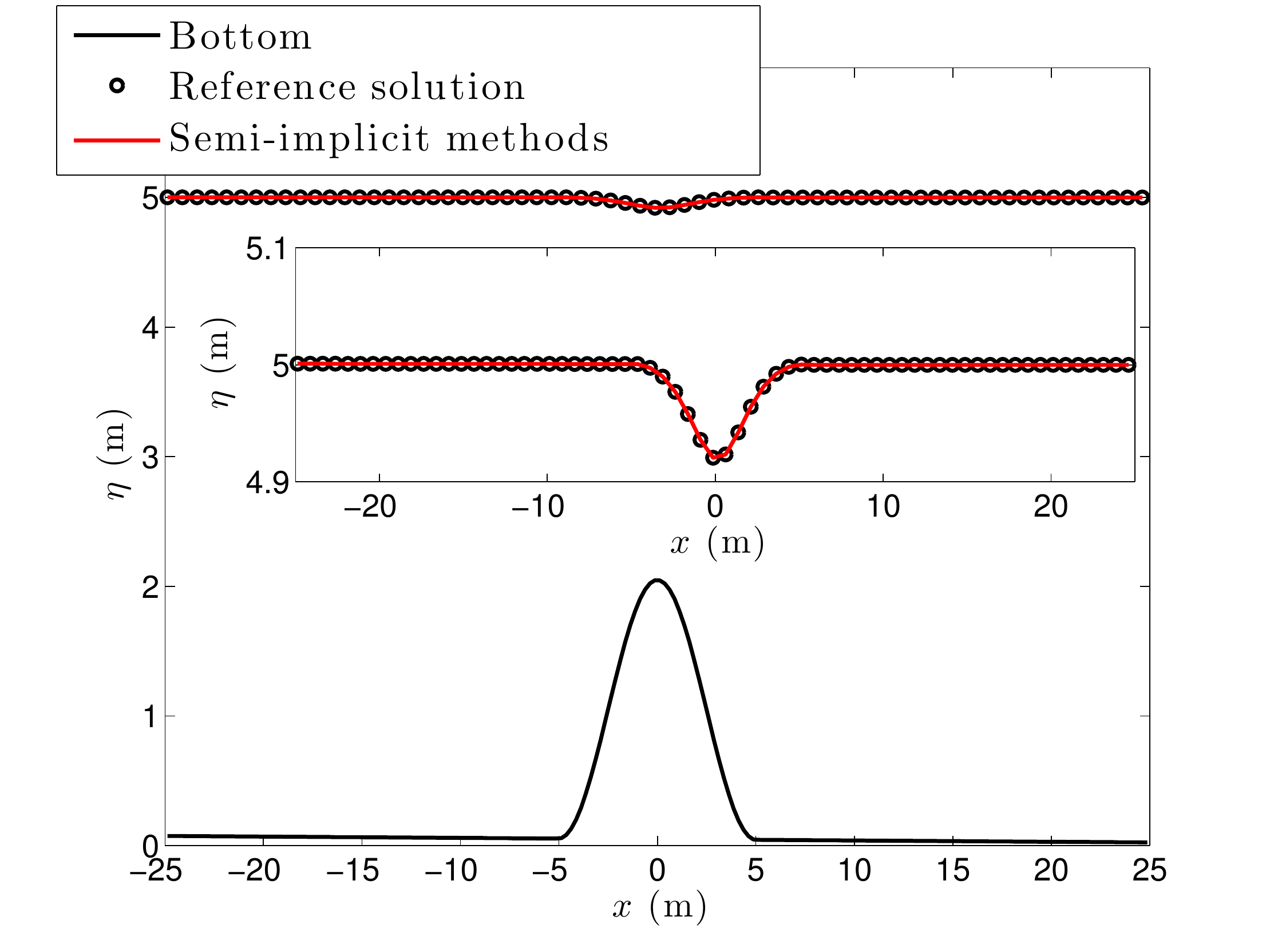}
  \caption{\label{fig:subc_freesurface} \footnotesize \textit{Free surface profile at steady state, as computed  in the steady subcritical flow test by the semi-implicit methods (solid red line) and reference solution (black circles) computed with the explicit scheme.  The inset figure is a zoom of the free surface profile.}}
 \end{figure}
 
\begin{table}
 \begin{center}
 \begin{tabular}{cccccccc}
\hline\\
SI-method & $\dt$ (s) & $C_{vel}$  & $C_{cel}$ & Err$_{\eta} $ [$l_2/l_\infty$]  & Err$_{u}$ [$l_2/l_\infty$] \\ 
 & & &  &   ${\small(\times10^{-6})}$ &  ${\small(\times10^{-5})}$\\ \hline\\
$\theta=0.55$ & 0.11 & 0.71 & 3.58 & 1.58/1.8 & 1.84/7.11\\
$\theta=0.7$ & 0.11 & 0.70 & 3.58 & 1.58/1.8  & 1.84/7.11\\
IMEX-ARK2 & 0.11 & 0.72 & 3.5 & 1.58/1.8  & 1.84/7.11\\
\hline\end{tabular}
 \caption{\footnotesize \it{Relative errors and Courant numbers achieved by using semi-implicit methods    in the steady subcritical flow test.}}
\label{tab:tabla2}
  \end{center}
 \end{table}
\begin{figure}[!h]\centering
\includegraphics[width=0.84\textwidth]{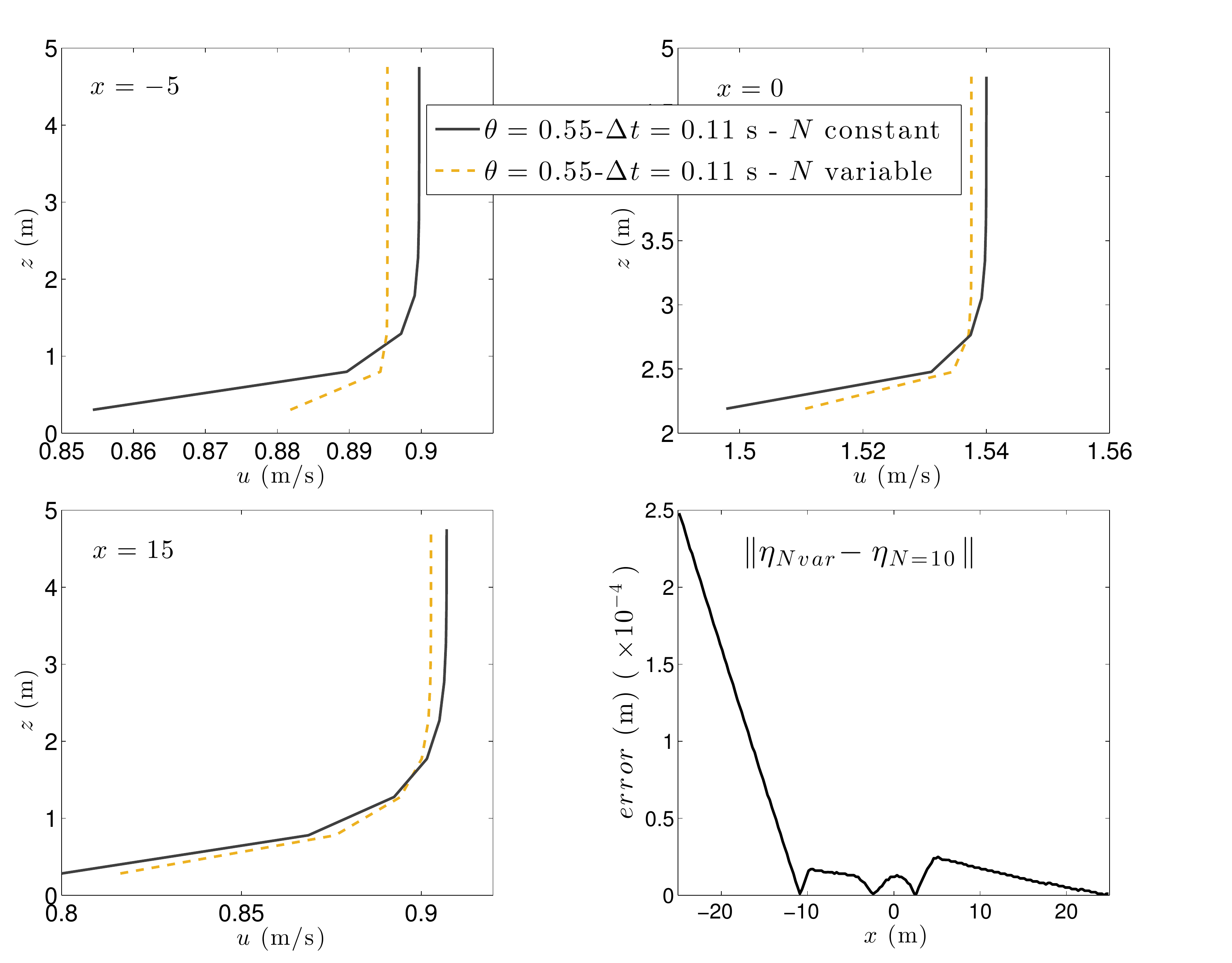}
  \caption{\label{fig:subc_nvar} \footnotesize \textit{Vertical profiles of horizontal velocity
  in the steady subcritical flow test, 
  obtained with the $\theta$-method  ($\theta = 0.55$) 
  and either   $10$ layers in the whole domain (solid black line) or a single layer in the first half of the domain  only (dashed yellow line). Profiles are taken at steady state at  the points $x = -5, 0, 15$ m. The solid black line denotes the absolute difference between the free surface computed with $10$ layers in the whole domain  or a single layer in the first half of the domain  only.}}
 \end{figure}

In figure \ref{fig:subc_freesurface} we see the free surface at the steady state, as computed with the semi-implicit $\theta$-method and IMEX-ARK2, along with the reference solution. In table \ref{tab:tabla2} we show the relative errors and the maximum Courant numbers achieved. The results computed with the semi-implicit methods are identical in this steady state case. Figure \ref{fig:subc_nvar} shows the absolute  difference on the free surface by using a semi-implicit method with either a constant number of layers or considering
\begin{equation}\label{eq:nvar2}
N = \left\{\begin{array}{ccl}
10 & &x> - 10,\\
1 & &\text{otherwise}.
\end{array}\right.
\end{equation}
The order of this difference is $10^{-4}$, with larger values where only one layer is employed. We also show the vertical profiles of horizontal velocity at three different points $x = -5, 0, 15$ m. These results show that we can reduce the number of degrees of freedom of our system  from 2210 to 1661, without a significant loss of accuracy where the multilayer configuration is kept.

\subsection{Tidal forcing over variable bathymetry}

\begin{figure}[!h]
\centering\includegraphics[width=0.45\textwidth]{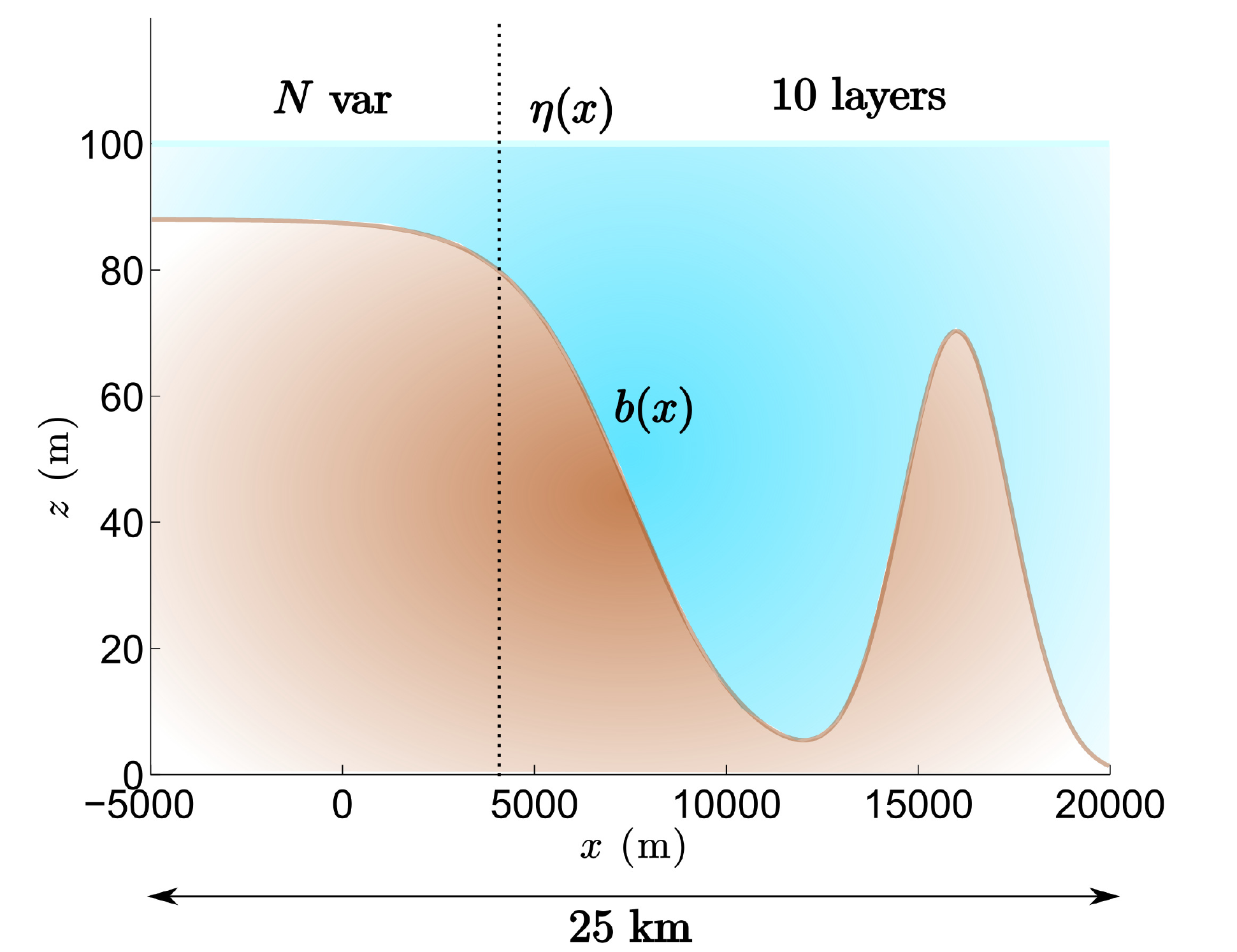}
  \caption{\label{fig:bat_marea} \footnotesize \textit{Sketch of the bottom topography.}}
 \end{figure}

 In this test we try to simulate a more realistic situation for coastal flow simulations. We consider a domain of length $L = 25$ km. The bottom bathymetry is taken as in Figure \ref{fig:bat_marea}, such that the  
 bathymetry is much shallower in one part of the computational
 domain than in the other. We define
 $$b(x) = z_0 - z_1\,\text{tanh}(\lambda\,(x-x_0)) +  70\,e^{-(x-x_1)^2/\sigma^2},$$
 with $z_0  = -z_1 = 44$, $x_0 = 7500$, $x_1 = 16000$, $\lambda = -1/3000$ and $\sigma=2000$. We
 consider   water at rest and constant free surface $\eta_0(x) = 100$ m at initial time. Subcritical boundary condition are imposed. The upstream condition is $q(-5000,t) = 1$ $m^2\,s^{-1}$, and the tidal downstream condition is $\eta(L,t) = 100 + 3\sin(\omega t)$ m, where $\omega = 2\pi/43200$. We simulate three 12-hours periods of tide, i.e., 36 hours. The friction parameters are taken as in previous tests with the exception of $\Delta z_0 = 3.3 \times10^{-3}$, which increases the bottom friction in order to obtain a more complex velocity field. In this case, a wind stress is included with a wind velocity of $1$ $m\,s^{-1}$. As in previous tests, we use $10$ vertical layers in the multilayer system and a uniform space discretization step $\dx = 50$ m.
 
\begin{figure}[h]
\centering\includegraphics[width=1\textwidth]{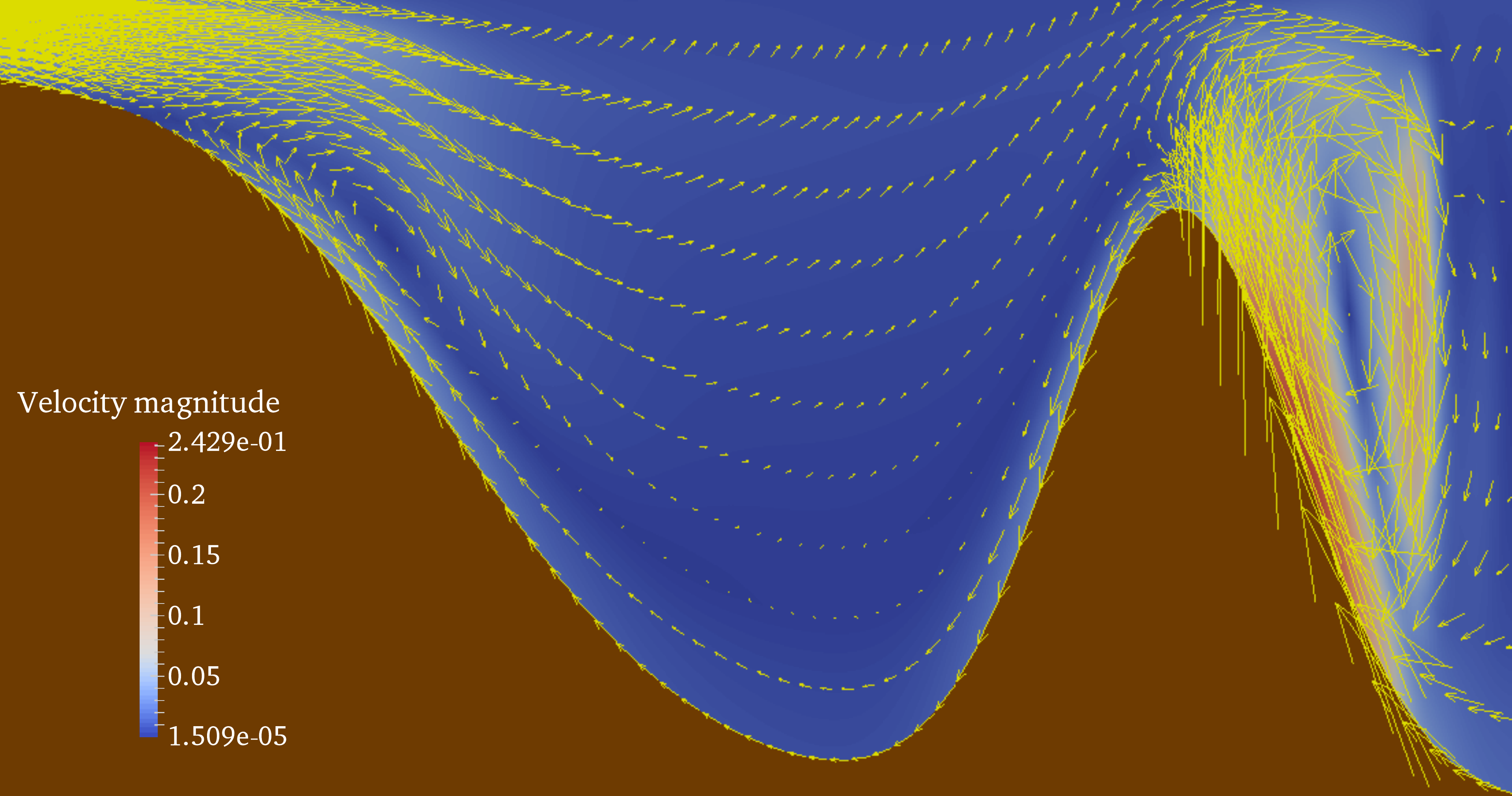}
  \caption{\label{fig:marea_vectorial} \footnotesize \textit{Vector map of the whole velocity field $\vu = (u,w)$ at time t = 33 h. Colors represent the magnitude of the velocity.}}
 \end{figure}

Figure \ref{fig:marea_vectorial} shows the obtained velocity field, where we can see some recirculations. Moreover, regarding the deepest area we realise that the upper and lower velocities has opposite direction.
   
\begin{figure}[!ht]
 \hspace{-1.4cm} \includegraphics[width=1.2\textwidth]{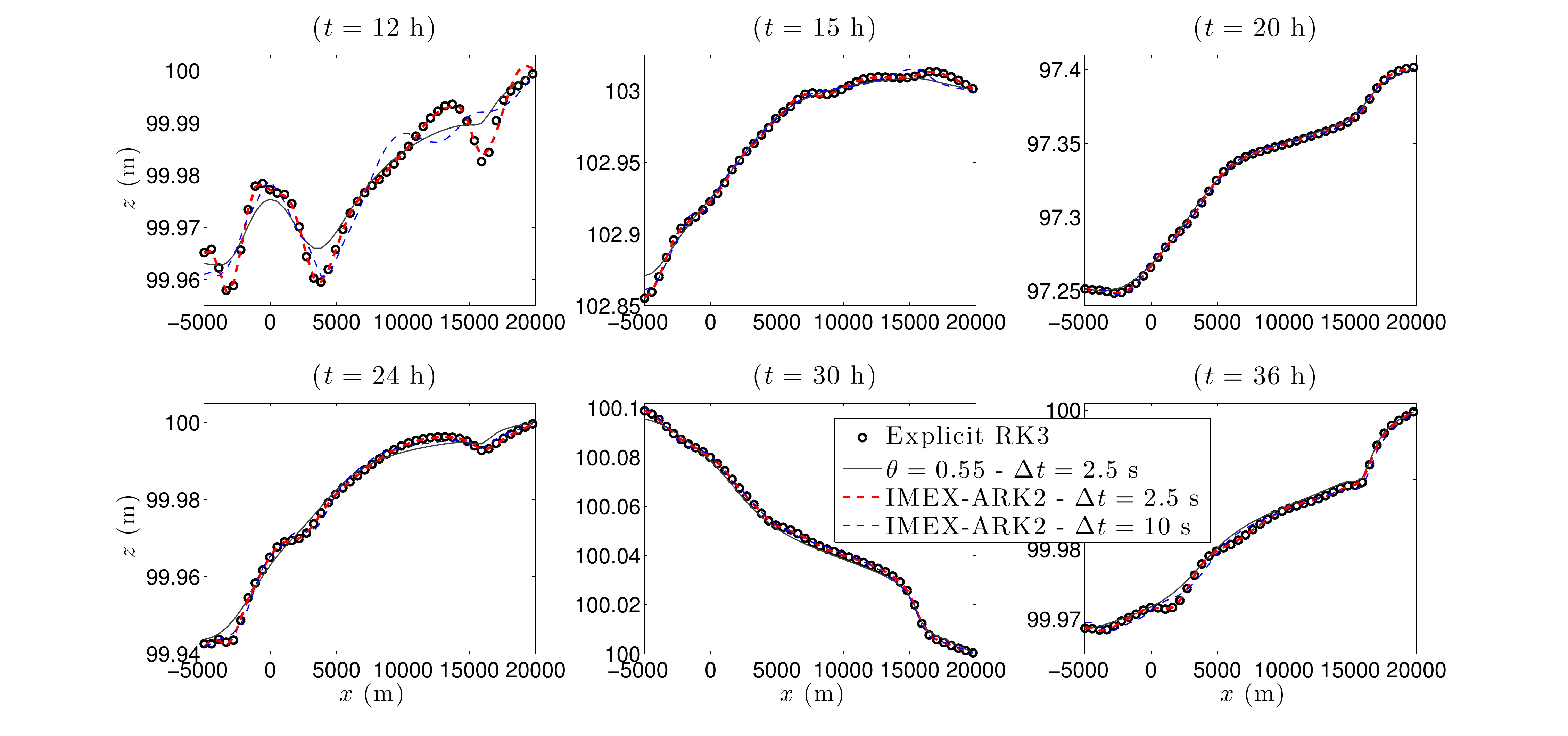}
  \caption{\label{fig:marea_eta} \footnotesize \textit{Free surface profile at different times 
  by using the semi-implicit methods (color lines) and the reference solution (black circles) computed with the explicit scheme in the tidal forcing test.}}
 \end{figure}

Figure \ref{fig:marea_eta} shows the free surface position at different times. We see that both  the $\theta$-method and the IMEX-ARK2 method are close of the reference solution. As in the free oscillation test, the IMEX-ARK2 approximates better the shape of the free surface. In particular, looking at table \ref{tab:tabla3}, where we report the relative errors at final time ($t =36$ h), we see that the second order method notably improves the results of the $\theta$-method. Note also that, 

\begin{table}[!ht]
\begin{center}
\begin{tabular}{cccccccc}
\hline\\
SI-method & $\dt$ (s) & $C_{vel}$  & $C_{cel}$ & Err$_{\eta} $ [$l_2/l_\infty$]  & Err$_{u}$ [$l_2/l_\infty$] \\ 
 & & & & ${\small(\times10^{-5})}$ & ${\small(\times10^{-2})}$\\ \hline\\
$\theta=0.55$ & 2.5 & 0.03 & 1.6 & 0.77/2.08  & 0.55/1.01\\
IMEX-ARK2 & 2.5 & 0.03 & 1.6 & 0.10/0.26  & 0.05/0.06\\
$\theta=0.55$ & 5 & 0.05 & 3.2 & 1.32/2.95  & 0.89/1.35\\
IMEX-ARK2 & 5 & 0.05 & 3.2 & 0.24/0.75 & 0.16/0.19\\
$\theta=0.55$ & 10 & 0.1 & 6.3 & 2.41/4.45 & 1.51/1.86\\
IMEX-ARK2 & 10 & 0.1 & 6.3 & 0.69/1.42  & 0.32/0.65\\
$\theta=0.55$ & 25 & 0.24 & 15.8 & 5.34/8.36 & 3.08/3.53\\
IMEX-ARK2 & 25 & 0.25 & 15.8 & 1.02/2.31  & 0.44/0.90\\
$\theta=0.55$ & 55 & 0.52 & 34.8 & 10.2/14.7 & 5.26/5.81\\
IMEX-ARK2 & 55 & 0.55 & 34.8 & 1.43/3.29  & 0.67/0.89\\
\hline\end{tabular}
\caption{\footnotesize \it{Relative errors and Courant numbers achieved by using semi-implicit methods at $t= 36$ h in the tidal forcing test.}}
\label{tab:tabla3}
\end{center}
\end{table}
\noindent in this typical coastal subcritical regime, large values of the Courant number can be achieved,   the maximum value being $C_{cel} = 34.8, $ without sensibly degrading the accuracy of the results. 
 
In table \ref{tab:tabla4} we report the computational times and speed-up achieved. With the explicit code  
about $16$ minutes of computation are required ($2.5$ hours for the reference solution), while the semi-implicit methods can reduce this time to seconds. Note also that  the IMEX-ARK2 is sensibly more efficient than the $\theta$-method in this case, since it is about $2.3$ times more expensive than the $\theta$-method, whereas the errors decrease by a much bigger factor. 
  
  \begin{table}[h]
 \begin{center}
 \begin{tabular}{cccccc}
 \hline\\
Method & $\dt$ (s) &  $C_{cel}$ & Comput. time (s) & Speed$-$up\\\hline\\
 Runge-Kutta 3 & -  & 0.1 (ref. sol.) & 9040 (150.6 m)  & -\\
 Runge-Kutta 3 & -  & 0.88 & 1014 (16.9 m) & 1\\
$\theta=0.55$ & 2.5 & 1.6 & 230 (3.8 m)  & 4.4\\
IMEX-ARK2 & 2.5 & 1.6 & 544 (9.1 m)  & 1.9\\
$\theta=0.55$ & 5 & 3.2 & 116 (1.9 m)  & 8.7\\
IMEX-ARK2 & 5 & 3.2 & 271 (4.5 m) & 3.74\\
$\theta=0.55$ & 10 & 6.3 & 58 & 17.5\\
IMEX-ARK2 & 10 & 6.3 &  136 (2.3 m) & 7.5\\
$\theta=0.55$ & 25 & 15.8 & 23 & 44.1\\
IMEX-ARK2 & 25 & 15.8 & 54  & 18,7\\
$\theta=0.55$ & 55 & 34.8 & 10 & 101.4\\
IMEX-ARK2 & 55 & 34.8 & 24 & 42.3\\
\hline\end{tabular}
 \caption{\footnotesize \it{Computational times and speed-up in the tidal forcing test case for the simulation up to $t=36$ h.}}
\label{tab:tabla4}
  \end{center}
 \end{table}

We also investigate the influence of simplifying the vertical discretization in the shallowest part of the domain (see figure \ref{fig:bat_marea}). We consider three different configurations, which we denote hereinafter  as (NVAR1)-(NVAR3). Firstly, we totally remove the vertical discretization by considering a single layer in the first part of the domain:
\begin{equation}\label{eq:marea_nvar1}
\tag{NVAR1}
N = \left\{\begin{array}{lll}	
10, & l_i = 1/10, \,i=1,...,N, & \mbox{ if } x\leq 4000;\\
1, & l_1 = 1,  & \mbox{ otherwise}.\\
\end{array}\right.
\end{equation}

Secondly, we keep a thin layer close to the bottom in order to improve the approximation of the friction term:
\begin{equation}\label{eq:marea_nvar2}
\tag{NVAR2}
N = \left\{\begin{array}{lll}	
10, & l_i = 1/10, \,i=1,...,N, & \mbox{ if } x\leq 4000;\\
2, & l_1 = 0.1, l_2 = 0.9  & \mbox{ otherwise}.\\
\end{array}\right.
\end{equation}
Finally, we improve again the vertical discretization close to the bottom by adding another thin layer:
\begin{equation}\label{eq:marea_nvar3}
\tag{NVAR3}
N = \left\{\begin{array}{lll}	
10, & l_i = 1/10, \,i=1,...,N, & \mbox{ if } x\leq 4000;\\
3, & l_1 = l_2 = 0.1, l_3 = 0.8,  & \mbox{ otherwise}.\\
\end{array}\right.
\end{equation}
 
\begin{figure}[!h]
 \hspace{-1.4cm} \includegraphics[width=1.2\textwidth]{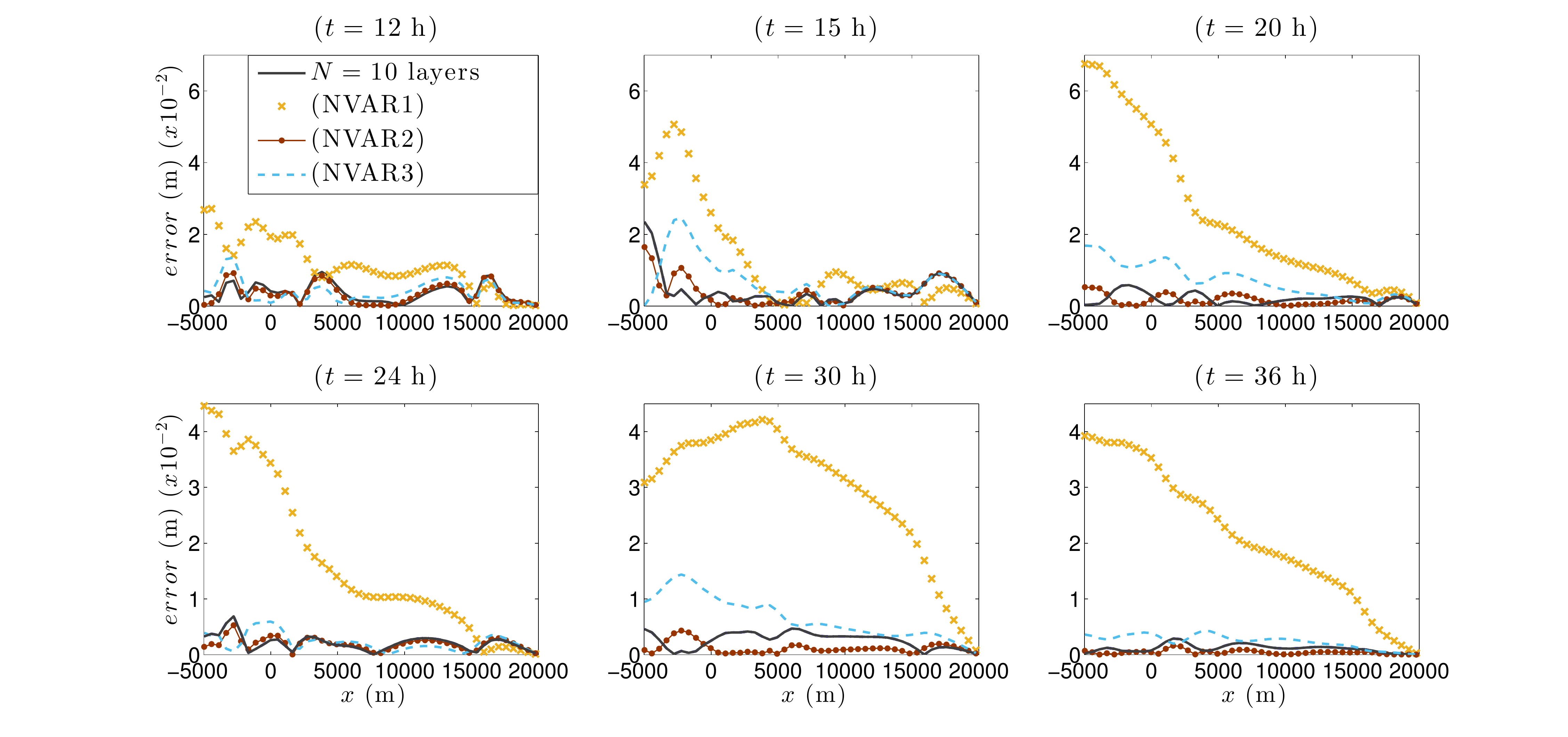}
  \caption{\label{fig:marea_nvar} \footnotesize \textit{Absolute errors for the free surface  at different times obtained in the tidal forcing test with the $\theta$-method ($\theta = 0.55$ and $\dt = 5$ s) and either $10$ layers in the whole domain (solid black line) or configurations (NVAR1)-(NVAR3) in the first part of the domain.}}
 \end{figure}

\begin{figure}[!h]
 \hspace{-1.4cm} \includegraphics[width=1.2\textwidth]{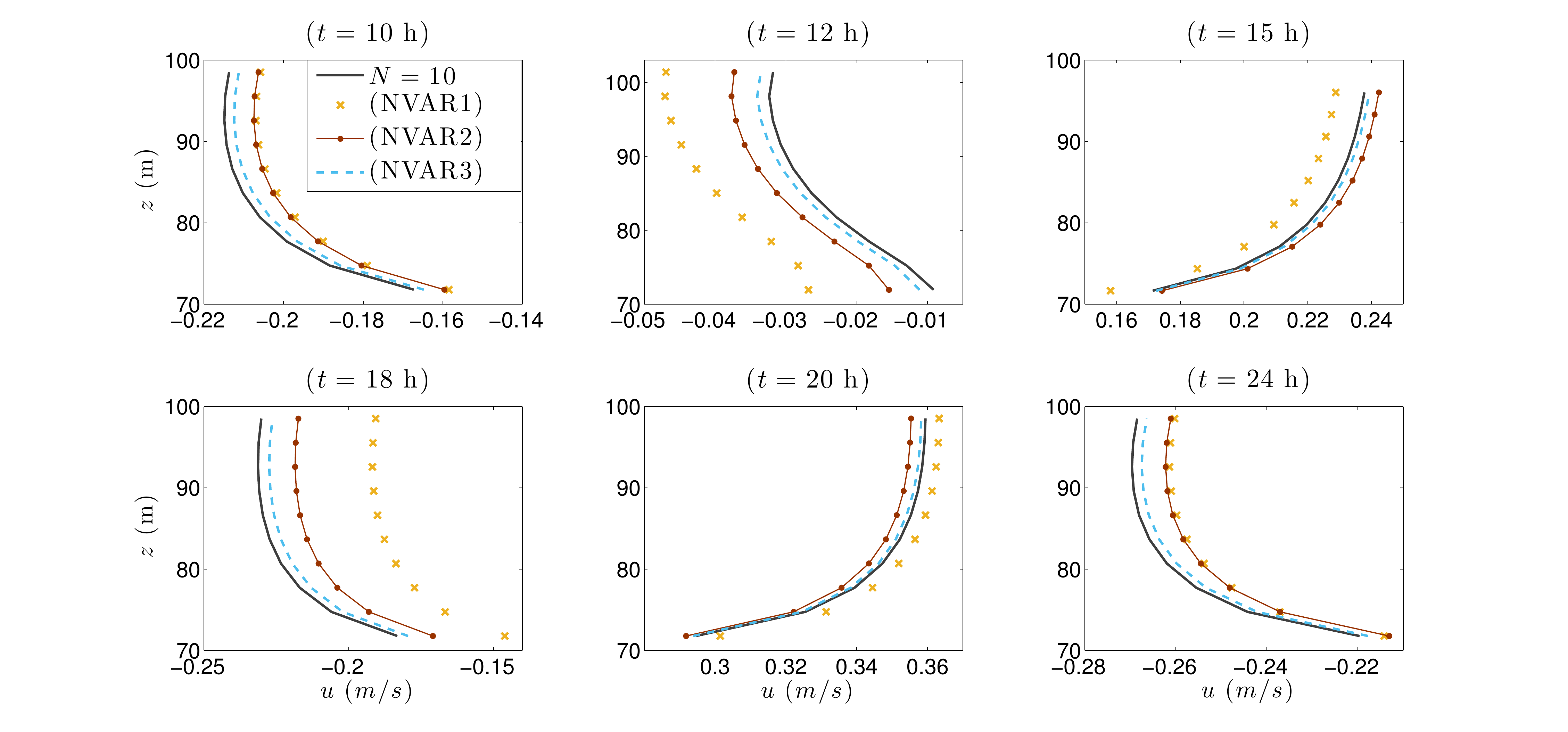}
  \caption{\label{fig:marea_nvar_vel} \footnotesize \textit{Vertical profiles of horizontal velocity obtained in the tidal forcing test with the $\theta$-method  ($\theta = 0.55$ and $\dt = 5$ s) and either $10$ layers in the whole domain (solid black line) or configurations (NVAR1)-(NVAR3) in the first part of the domain. Profiles are taken at the point $x = 16025$ m and times $t = 12,\,15,\,20,\,24,\,30,\,36$ h.}}
 \end{figure}

In this way, the number of degrees of freedom of the multilayer system is reduced from 5510 to 3890 (NVAR1),  4070 (NVAR2), or 4250 (NVAR3). Note that configurations (NVAR2) and (NVAR3) employ a non-uniform distribution of the vertical layers. Figure \ref{fig:marea_nvar} shows the absolute errors with the $\theta$-method with $\dt=5$ s ($C_{cel} =3.2$) using $10$ layers in the whole domain and with configurations (NVAR1)-(NVAR3). We see that  the simplest configuration (NVAR1) leads to the largest error. However, by using configurations (NVAR2) and (NVAR3) these errors are much more similar to the case in which a constant number of layer is employed in the whole domain. As expected, the  smallest error is achieved with the configuration (NVAR3). Figure \ref{fig:marea_nvar_vel} shows the vertical profile of horizontal velocity at point $x=16025$ m (the top of the peak) at different times. The conclusions are similar, i.e., the differences are larger with  configuration (NVAR1), whereas (NVAR2) and (NVAR3) give accurate approximations of the vertical profile obtained with a constant number of layers.

 \subsection{An application to sediment transport problems}
 
 In order to emphasize the usefulness of the proposed method and the potential advantages
 of its application   to more realistic problems, we consider the extension of equations  \eqref{eq:sistema3D_u}
to the  movable bed case. For simplicity, we work with a decoupled, essentially monophase model,
according to the classification in \cite{garegnani:2011}, \cite{garegnani:2013}, which is appropriate
 in the limit of small sediment concentration. Quantity $b $ in \eqref{eq:sistema3D_u} is then assumed to be dependent on time and an Exner equation for the bed evolution is also considered
 \begin{equation}\label{eq:exner}
\frac{\partial b}{\partial t}+  \frac{\partial Q_b}{\partial t}   = 0,
\end{equation}
where $\xi = 1/(1-\rho_0)$ with $\rho_0$ the porosity of the sediment bed, and the solid transport discharge is defined by an appropriate formula, see e.g. \cite{castro:2008}. Here we consider a simple definition of the solid transport discharge given by the Grass equation $$Q_b = A_g u^3,$$ where $A_g \in (0,1)$ is an experimental constant depending on the grain diameter and the kinematic viscosity. 
For control volume $i$, equation \eqref{eq:exner} is easy discretized along the lines of section \ref{si_td}. For the $\theta$-method, the discrete equation reads 
\begin{eqnarray}
z_{b,i}^{n+1} = z_{b,i}^n &+&\theta\,\xi\, A_g\,\dfrac{\dt}{\dx}\left( |u_{1,i-\frac12}^{n+1}|^2\,u_{1,i-\frac12}^{n+1} - |u_{1,i+\frac12}^{n+1}|^2\,u_{1,i+\frac12}^{n+1} \right)\\
&+&(1-\theta)\,\xi\, A_g\,\dfrac{\dt}{\dx}\left( |u_{1,i-\frac12}^{n}|^2\,u_{1,i-\frac12}^{n} - |u_{1,i+\frac12}^{n}|^2\,u_{1,i+\frac12}^{n} \right).\nonumber
\end{eqnarray}

On the other hand, the IMEX-ARK2 discretization of equation \eqref{eq:exner} consists of a simple updating of the values of the movable bed, since the values $u_{\a}^{n,j}$ are known when $z_b^{n,j}$ is computed. For the first stage we have $z_{b,i}^{n,1} =z_{b,i}^{n}$. Next, $z_{b,i}^{n,2}$ and $z_{b.i}^{n,3}$ are computed by the formula
\begin{eqnarray}
z_{b,i}^{n,j} = z_{b,i}^n &+&   \xi\, A_g\,\dfrac{\dt}{\dx}\,\dsum_{k=1}^{j}\tilde{a}_{jk}\,\left( | u_{1,i-\frac12}^{n,k}|^2\,u_{1,i-\frac12}^{n,k} - |u_{1,i+\frac12}^{n,k}|^2\,u_{1,i+\frac12}^{n,k} \right).\nonumber
\end{eqnarray}
Finally, the solution at time $n+1$ is
\begin{eqnarray}
z_{b,i}^{n+1} = z_{b,i}^n &+&   \xi\, A_g\,\dfrac{\dt}{\dx}\dsum_{j=1}^{3}\tilde{b}_{j}\,\left( |u_{1,i-\frac12}^{n,j}|^2\,u_{1,i-\frac12}^{n,j} - |u_{1,i+\frac12}^{n,j}|^2\,u_{1,i+\frac12}^{n,j} \right).\nonumber
\end{eqnarray}

We consider a simple test in which a parabolic dune is displaced by the flow  (see \cite{castro:2008}). The computational domain has length $1000$ m and $150$ nodes are used in the spatial discretization. We set the constant $A_g$ in the Grass formula to $0.001$ and we take the porosity value  $\r_0 = 0.4$. We consider viscosity effects with the same parameters as in the previous tests, disregarding wind stress. Subcritical boundary condition are imposed. The upstream condition is $q(0,t) = q_0(x)$ and the downstream one is $\eta(L,t) = 15$ m. The initial condition for the bottom profile is
given by
\begin{equation}\label{eq:dune}
z_{b,0} (x) = \left\{\begin{array}{ll}
0.1 + \sin^2\left(\dfrac{\pi(x-300)}{200}\right) & \mbox{if } 300\leq x \leq 500;\\
0.1 & \mbox{otherwise}, 
\end{array}\right.
\end{equation}
and the initial height is $h_0(x) = 15 - z_{b,0}(x)$. For the discharge, we take into account the vertical structure of the flow in order to have a single dune moving along the domain. With this purpose, we run a first simulation  of the movement of the dune \eqref{eq:dune}, where the initial discharge is  $q_i=15\ m^2\,s^{-1}$, for $i=1,\dots,N$, until it reaches a steady structure at the outlet. These values of the discharge are used as initial and upstream boundary condition in the final simulation. If a constant discharge were used, this would  sweep along the sediment in the initial part of the domain and create another dune within the computational domain. While this is physically correct, we prefer in this test to study a simpler configuration.

\begin{figure}[!h]
\centering\includegraphics[width=0.8\textwidth]{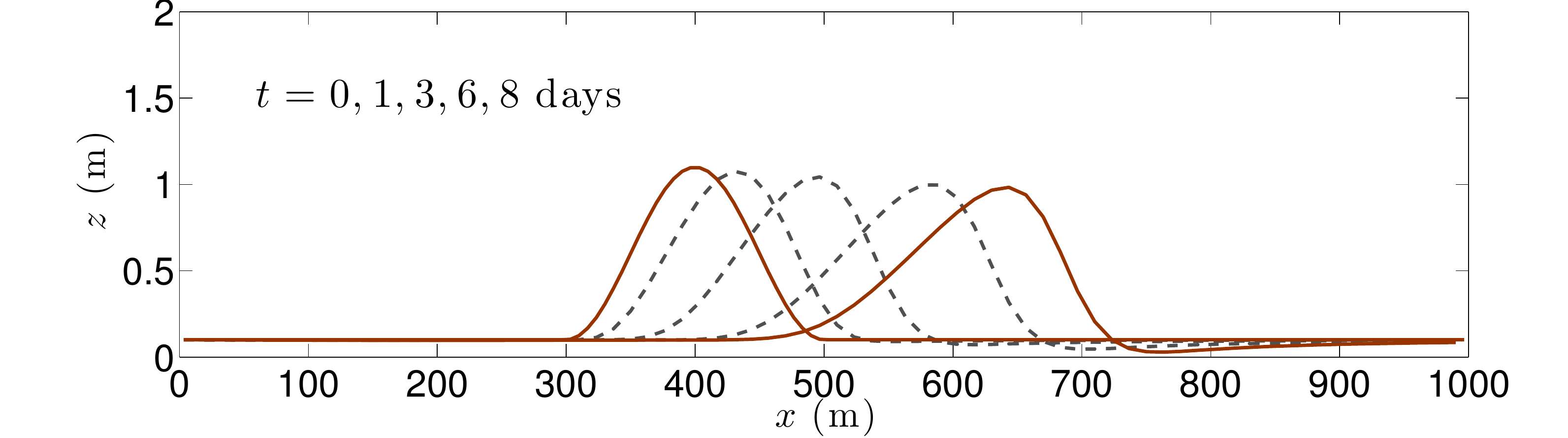}
  \caption{\label{fig:sedi_evol} \footnotesize \textit{Profile of the dune at different times
  in the sediment transport test case, including the initial condition and the final position.}}
 \end{figure}

\begin{figure}[!h]
 \hspace{-1.4cm}\includegraphics[width=1.2\textwidth]{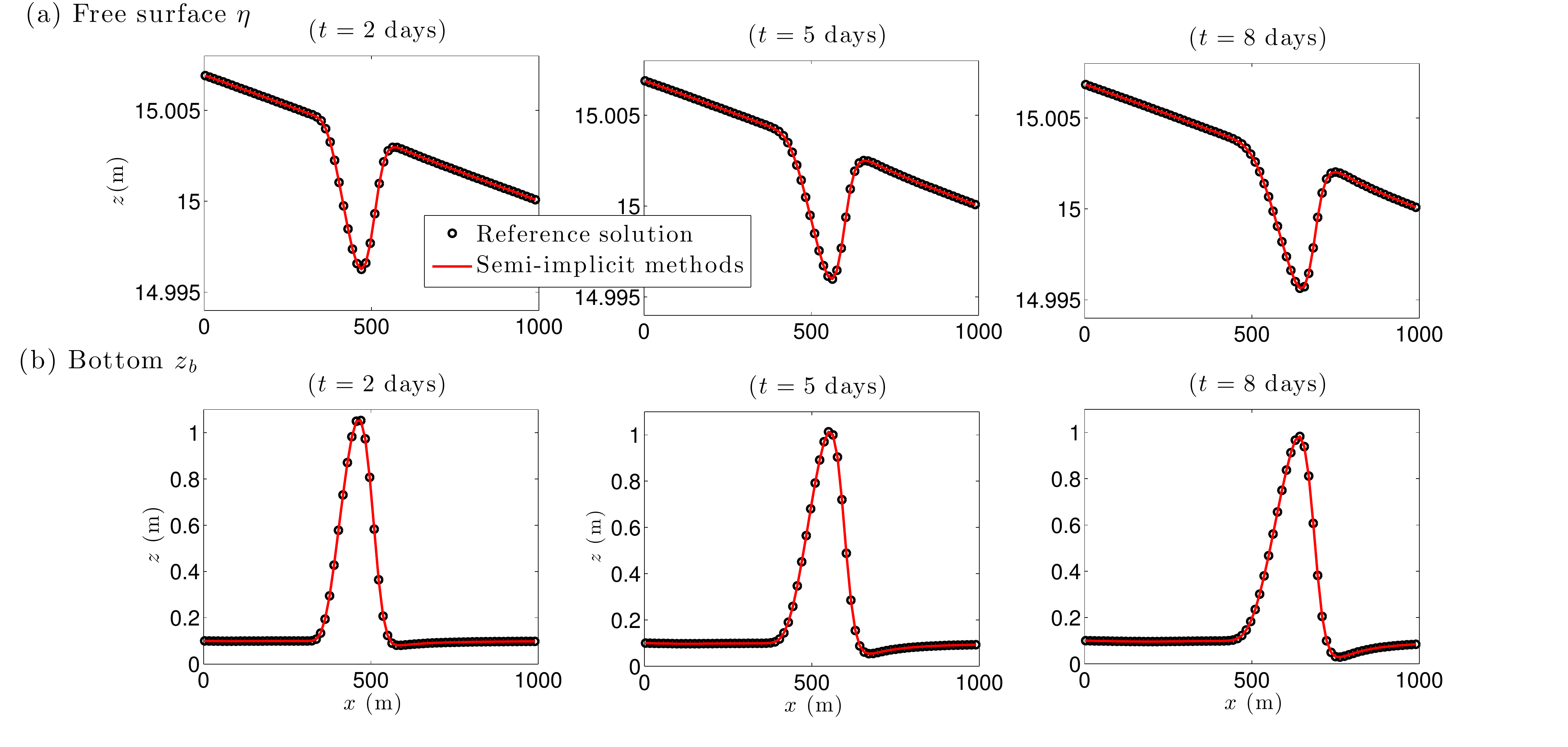}
  \caption{\label{fig:sedi_eta_zb} \footnotesize \textit{Free surface and bottom profile at different times in the sediment transport test case, as computed by  the semi-implicit methods (solid red line) and by the reference  explicit scheme (black circles).}}
 \end{figure}

We use 10 layers in the multilayer code and simulate until $ t= 691200$ s ($8$ days). Figure \ref{fig:sedi_evol} shows the evolution of the dune and figure \ref{fig:sedi_eta_zb} shows zooms of evolution of the free surface and of the movable bed, as computed with either the explicit third order Runge-Kutta or the  semi-implicit ($\theta$-method and IMEX-ARK2). The results are essentially indistinguishable. This is confirmed looking at table \ref{tab:tabla5}, where we report the relative errors and the Courant number achieved. We see that there are not significant differences between the semi-implicit methods, due to the fact that the flow is essentially a steady
one and the bed evolution is very slow.

 
  \begin{table}[h!]
 \begin{center}
 \begin{tabular}{cccccccc}
 \hline\\
SI-method & $\dt$ (s) & $C_{vel}$ & $C_{cel}$ & Err$_{\eta} $ [$l_2/l_\infty$]  & Err$_{u}$ [$l_2/l_\infty$] & Err$_{b} $ [$l_2/l_\infty$] \\ 
&  & &  & ${\small(\times10^{-7})}$ &  ${\small(\times10^{-6})}$  & ${\small(\times10^{-5})}$\\ \hline\\
$\theta=0.55$ & 1 & 0.17 & 1.98 & 1.4/5.35  & 0.29/1.41 & 1.09/1.52\\
IMEX-ARK2 & 1 & 0.16 & 1.97 & 1.29/5.39  & 0.27/1.41 & 1.03/1.40\\
$\theta=0.55$ & 2 & 0.34 & 3.94 & 1.69/6.13  & 0.55/2.90 & 2.25/3.13\\
$\theta=0.6$ & 2 & 0.34 & 3.94 & 1.69/6.13  & 0.55/2.90 & 2.25/3.13\\
IMEX-ARK2 & 2 & 0.33 & 3.93 & 1.68/6.47 & 0.50/2.33 & 2.11/2.87\\
\hline\end{tabular}
 \caption{\footnotesize \it{Relative errors and Courant numbers achieved in the sediment transport test case by   semi-implicit methods at $t= 192$ hours (eight days).}
  }
\label{tab:tabla5}
  \end{center}
  
 \end{table}
 
 As remarked before, a rigorous comparison of the efficiency of the proposed methods is not possible in our preliminary implementation. However, a preliminary assessment is reported in Table \ref{tab:tabla6}, showing the computational time and the speed-up obtained for the simulation of 192 hours (8 days). For the reference solution with the explicit scheme approximately  13 hours are necessary (78 minutes with maximum $C_{cel}$), whereas 8 minutes (respectively, 19 minutes) are needed with the $\theta$-method and IMEX-ARK2 method when considering a time step $\dt =2$ s. This gives a speed up of 9 (4 for the IMEX-ARK2). Even taking a small time step ($\dt=1$ s) the computational time required is notably reduced to 17 min (39 min for the IMEX).

  \begin{table}
 \begin{center}
 \begin{tabular}{cccccc}
 \hline\\
Method & $\dt$ (s) &  $C_{cel}$ & Comput. time (s) & Speed$-$up\\\hline\\
 Runge-Kutta 3 & -  & 0.1 (ref. sol.) & 45978 (12.7 h)  & -\\
 Runge-Kutta 3 & -  & 0.99 & 4700 (78.33 m)  & 1\\
$\theta$-method & 1  & 1.98 & 1048 (17.5 m) & 4.5 \\
IMEX-ARK2 & 1  & 1.97 & 2368 (39.4 m) & 1.99 \\
$\theta$-method & 2 & 3.94 & 509 (8.5 m) & 9.2\\
IMEX-ARK2 & 2 & 3.93 & 1164 (19.4 m) & 4.04 \\
\hline\end{tabular}
 \caption{\footnotesize \it{Computational times and speed-up in the sediment transport test case for the simulation up to $t=192$ h (eight days).}}
\label{tab:tabla6}
  \end{center}
 \end{table}

\begin{figure}[!h]
\centering\includegraphics[width=0.6\textwidth]{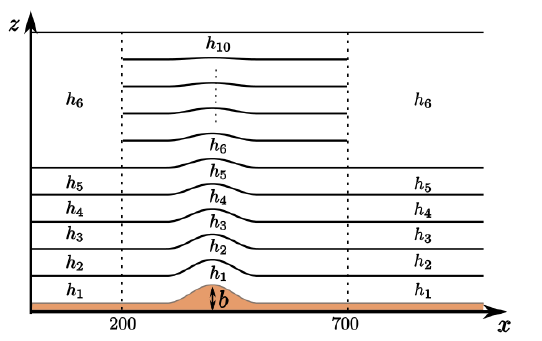}
  \caption{\label{fig:sedi_capas} \footnotesize \textit{Sketch of the multilayer configuration with the variable number of layers for the sediment transport test case.}}
 \end{figure}

 Finally, we can further reduce the computational time by reducing locally the number of layers employed. In this test, the vertical structure cannot be completely removed without causing a major loss of accuracy, since the dynamics of the movable bed depends on the velocity of the layer closest  to the bottom. For this reason, we consider the following configuration (see also figure \ref{fig:sedi_capas}): 

\begin{equation}\label{eq:nvar_sedi}
N = \left\{\begin{array}{lll}	
10, & l_i = 1/10, \,i=1,...,N, & \mbox{ if }200\leq x\leq 700;\\
6, & l_i = 1/10, \,i=1,...,5; l_6 = 0.5,  & \mbox{ otherwise}.\\
\end{array}\right.
\end{equation}

\begin{figure}[!h]
 \hspace{-1.4cm} \includegraphics[width=1.2\textwidth]{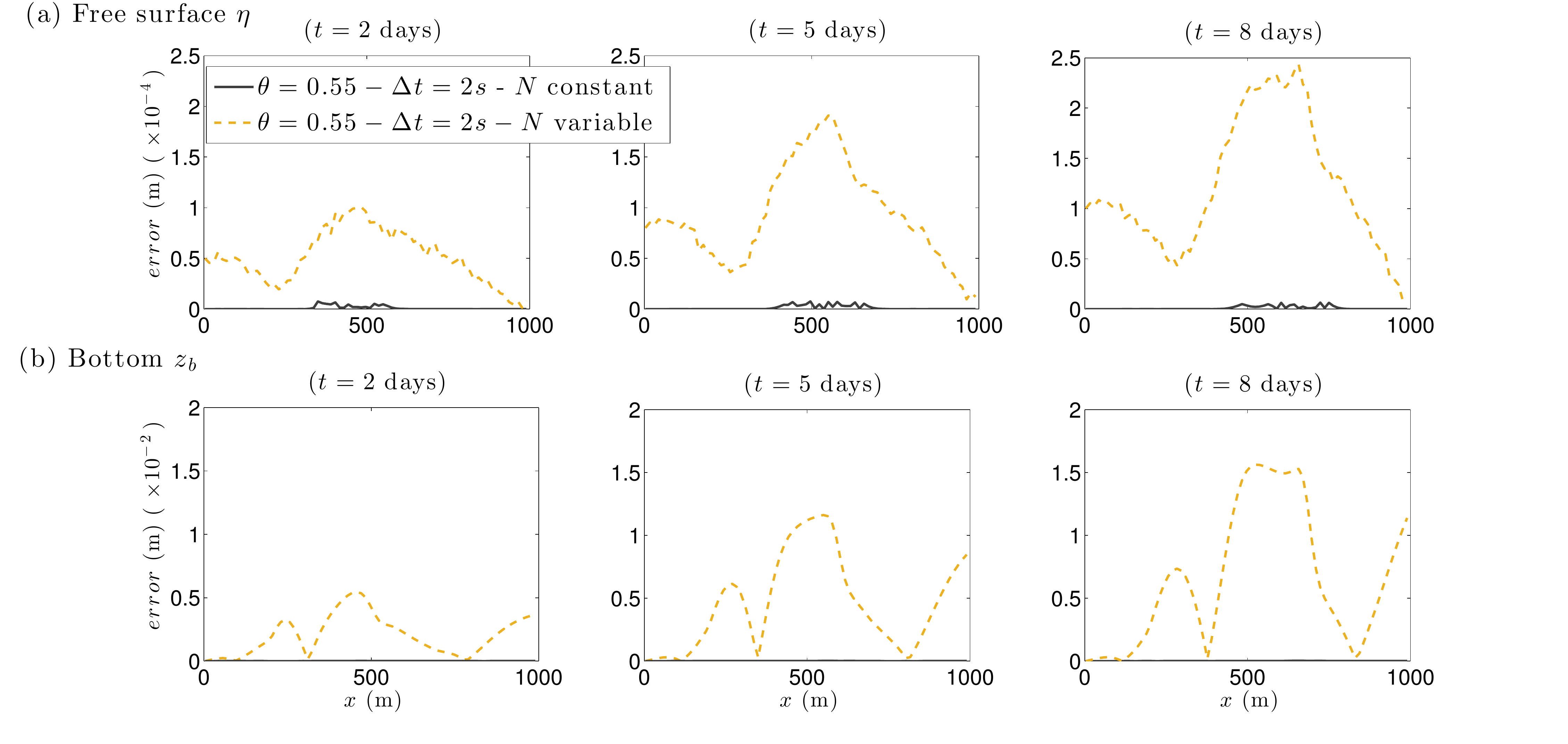}
  \caption{\label{fig:sedi_error_nvar} \footnotesize \textit{Absolute differences for the free surface ($\eta$) and bottom ($z_b$) at different times in the sediment transport test case, by using the $\theta$-method  ($\theta = 0.55$ and $\dt = 2$ s). We compare the results with $10$ layers in the whole domain (solid black line) with those obtained with the variable number of layers (see \eqref{eq:nvar_sedi}, dashed yellow line).}}
 \end{figure}

Note that, in this way, both a variable number of  vertical layers and a non-uniform distribution of these layers are tested. Figure \ref{fig:sedi_error_nvar} shows the absolute differences on the free surface and on the movable bed profiles at different times when we use wither a constant number of layers ($N=10$) or the configuration \eqref{eq:nvar_sedi}. The difference between both configurations for the bottom is lower than the 2\% of its thickness, whereas the number of degrees of freedom of the problem is reduced  from 1660 to 1352.

  \section{Conclusions}
 \label{conclu} \indent
 
We have proposed two concurrent strategies to make multilayer models more efficient and fully competitive with their $z-$ and $\sigma-$coordinates counterparts. On one hand, we have shown how the number of vertical layers that are employed  can be allowed to vary over the computational domain. Numerical experiments  show that, in the typical regimes in which the application of multilayer shallow water models is justified, the resulting discretization does not introduce any major spurious feature and allows to reduce substantially the computational cost in areas with complex bathymetry.  Furthermore, efficient semi-implicit discretizations have been applied for the first time to this kind of models, allowing to achieve significant computational gains in subcritical regimes. This makes multilayer discretizations fully competitive  with $z-$coordinate discretizations for large scale, hydrostatic flows.
In addition, a more efficient way to implement the IMEX-ARK method to discretize the multilayer system, which mimics what done for simpler $\theta$-method, has been proposed.
In particular, in the applications to tidally forced flow and to the sediment transport problem,  we have shown that the computational time required is significantly reduced and that the vertical number of layers, as well as their distribution, can be adapted to the local features of the problem.

In future work,  we will be interested in applying this approach to more realistic simulations.
In particular, we will extend the proposed approach to variable density flows in the Boussinesq regime.
Furthermore, we plan to couple multilayer vertical discretizations
  to the adaptive, high order horizontal discretizations proposed in 
\cite{tumolo:2015}, \cite{tumolo:2013}, in order to achieve maximum accuracy for the envisaged application regimes.

\section*{Acknowledgements}
This work was partially supported by the Spanish Government and FEDER through the research projects MTM2012-38383-C02-02 and  MTM2015-70490-C2-2-R. Part of this work was carried out during visits by J. Garres-D\'iaz at MOX Milano and L. Bonaventura at IMUS Sevilla. 
\bibliographystyle{plain}
\bibliography{multilayer}

\end{document}